\input amstex






\font\rm=cmr10 \rm

\font\bf=cmb10
\font\Rm=cmr9 at 11pt
\rm
\font\it=cmsl9 at 10pt
 at 7pt

\font\Rrm=cmr17 at 16pt
   \font\Rm=cmr12 at 11.5pt

\long\def\Pf{\par\noindent {\it Proof.} }
\def\({\left(}
\def\){\right)}
\def\st{such that }
\def\qed{\hfill$\bullet$\vskip 4pt}

\def\brcs#1{\left\{ #1\right\}}

\def\comp{\circ}
\def\Log{\text{Log\,}}
\def\iso{\cong}
\def\wrt{with respect to }
\def\:{\,:}

\def\supp{\text{supp}\,}

\def\det{\text{det\,}}

\def\ker{\text{ker\,}}

\def\C{\text{\bf C}}

\def\I{\text{I\,}}

\def\tr{\text{\,tr}}
\def\Tr{\text{\,Tr}}
\def\R{\text{\bf R}}
\def\N{\text{\bf N}}
\def\Z{\text{\bf Z}}

\def\Arrow #1;#2.{#1\:#2 \to }

\def\Set#1#2{\brcs{#1 \left|\vphantom{#1 #2} \right.#2}}


\def\oh#1{{\pmb o}\(#1\)}

 \def\supp#1{\text{supp\,}#1}
\def\Rrr#1,#2{{\Cal J}_{#1,#2}}
\def\slfrac#1#2{{\raise -.07 ex\hbox{$^{#1}$}}\!/\raise .35 ex \hbox{${}_{#2}$}}
\def\ssf #1/#2{\slfrac {#1}{#2}}

\def\pd #1,#2.{\frac {\partial #1}{\partial #2}}

   \long\def\Lem
#1.#2\par{\vskip4pt{\baselineskip=13pt\font\it=cmsl12 at
11.5pt\Rm
   \noindent {\rm \uppercase{#1}} #2\vskip3pt

   }} 

\long\def\Proclaim #1.#2 \endproclaim{\vskip4pt{\baselineskip=13pt\font\it=cmsl12 at
11.5pt\Rm
   \noindent {\rm \uppercase{#1}} #2\vskip3pt

   }} 

\long\def\remark #1\endremark{\vskip 2pt \noindent {\it Remark\/} #1\par}

\long\def\Sectionhead #1.#2:\par #3{\vskip 4pt \noindent {\bf #1 #2}vskip 2pt\noindent\nospace #3}

\long\def\Title #1\par {\noindent{\Rrm #1}\vskip 9pt}

 \long\def\SubT #1.{\noindent {\it #1\/} } 
 
 \long\def\SecT
#1\par{\vskip 3pt \noindent {\bf #1}\vglue1pt
   \noindent}

\long\def\subtitle #1.{\vskip 2pt \noindent {\it #1}}

\long\def\Rmk#1\par{\vskip 1pt \noindent {\it
Remark.} #1\vskip2pt}

\long\def\Abstract #1\par{{\leftskip= 3 true cm \rightskip = 3 true cm \font\it=cmsl10 \font\rm=cmr10 \baselineskip = 10pt
\parindent=.35 true cm\rm\noindent 
{\it Abstract} #1\vskip 8pt

}}

\long\def\Author #1 \par{\noindent{\it #1}}

 
\scrollmode\NoBlackBoxes
\magnification=1100
\long\def\Abstract #1\par%
{\vskip .2 true cm{\leftskip 1 true in \rightskip 1 true in \font\rm=cmr8 \rm
\baselineskip=1pt \font\it=cmsl8 \font\bf=cmb10 at 8pt
\parindent=0em {\bf Abstract} #1

}}
\comment
\font\rm=Times at 10pt

\font\bf=TimesB
\font\Rm=Times at 11pt
\rm
\font\it=TimesI at 10pt
\endcomment

\long\def\Pf{\par\noindent {\it Proof.} }
\def\({\left(}
\def\){\right)}
\def\st{such that }
\def\qed{\hfill$\bullet$\vskip 4pt}

\def\brcs#1{\left\{ #1\right\}}
\def\Set#1#2{\brcs{#1 \left|\vphantom{#1 #2} \right.#2}}

\def\C{\text{\bf C}}

\def\I{\text{I\,}}

\def\tr{\text{\,tr}}
\def\Tr{\text{\,Tr}}

\def\comp{\circ}
\def\Log{\text{Log\,}}
\def\iso{\cong}
\def\wrt{with respect to }
\def\:{\,:}
\def\Arrow #1;#2.{#1\:#2 \to }

\def\oh#1{{\pmb o}\(#1\)}
\def\R{\text{\bf R}}
\def\N{\text{\bf N}}
\def\Z{\text{\bf Z}}

 \def\supp#1{\text{supp\,}#1}
\def\Rrr#1,#2{{\Cal J}_{#1,#2}}

\def\slfrac#1#2{{\raise -.07 ex\hbox{$^{#1}$}}\!/\raise .35 ex \hbox{${}_{#2}$}}
\def\ssf #1/#2{\slfrac {#1}{#2}}

\def\pd #1,#2.{\frac {\partial #1}{\partial #2}}


   \long\def\Title #1\par {\noindent{\Rrm #1}\vskip 9pt}
 \long\def\SubT #1.{\noindent {\it #1\/} }   \long\def\SecT
#1\par{\vskip 3pt \noindent {\bf #1}\vglue1pt
   \noindent}
\long\def\subtitle #1.{\vskip 2pt \noindent {\it #1}}

\long\def\Rmk#1\par{\vskip 1pt \noindent {\it
Remark.} #1\vskip2pt}


\def\oneone{\One.1}
\def\onetwo{\One.2}
\def\onethr{1.3}
\def\onefou{1.4}
\def\onefiv{1.5}

\def\twoone{3.1}
\def\twotwo{3.2}
\def\twothr{3.3}
\def\twofou{3.4}
\def\twofiv{3.5}
\def\twosix{3.6}

\def\twoeig{3.8}
\def\twonin{3.9}

\def\throne{2.1}

\def\foutwo{4.2}
\def\fouthr{4.3}

\def\sixone{6.1}
\def\sixtwo{6.2}
\def\sixthr{6.3}

\def\sevthr{7.3}
\def\sevfou{7.4}
\def\sevfiv{7.5}
\def\sevsix{7.6}


\input diagrams.tex

\def\oneone{1.1}
\def\onetwo{1.2}
\def\onethr{1.3}
\def\onefou{1.4}
\def\onefiv{1.5}

\def\twoone{2.1}
\def\twotwo{2.2}
\def\twothr{2.3}
\def\twofou{2.4}
\def\twofiv{2.5}
\def\twosix{2.6}

\def\twoeig{2.8}
\def\twonin{2.9}
\def\twoten{2.10}

\def\throne{3.1}

\def\foutwo{4.2}
\def\fouthr{4.3}

\def\sixone{6.1}
\def\sixtwo{7.2}
\def\sixthr{7.1}

\def\sevthr{7.3}
\def\sevfou{7.4}
\def\sevfiv{7.5}
\def\sevsix{7.6}

\def\I{\text{I\,}}
\def\tr{\text{tr\,}}

\def\flo #1{\lfloor #1 \rfloor}

\let\iso=\cong

\def\thmxx{fill me in!!}

\def\tripnorm #1xxx{\left\|\hskip-.2em{\left| #1 \right|}\hskip-.2em\right\|}

\def\flo #1{\lfloor #1 \rfloor}

\def\diag{\text{diag}\,}

\Title Nearly approximate transitivity (AT) for circulant matrices%

\Abstract By previous work of Giordano and the author, ergodic actions of $\Z$ (and other discrete groups) are completely classified measure-theoretically by their dimension space, a construction analogous to the dimension group used in C*-algebras and topological dynamics. Here we investigate how far from AT (approximately transitive) can actions be which derive from circulant (and related) matrices. It turns out not very: although non-AT actions can arise from this method of construction, under very modest additional conditions, ATness arises; in addition, if we drop the positivity requirement in the isomorphism of dimension spaces, then all these ergodic actions satisfy an analogue of AT. Many examples are provided.

\Author David Handelman %
\plainfootnote{$^{1}$}{\rm Supported in part by an NSERC Discovery Grant.}

\SecT Introduction

Let $(X,\mu,T) $ (often abbreviated $T$) be a measure space with an invertible measurable transformation $T$ (up to sets of measure zero), \st $\mu \comp T^k$ is absolutely continuous \wrt $\mu$ for all integers $k$; we also assume the action of $T$ is ergodic, that is, any $T$-invariant set has measure zero or one.

Motivated by the fundamental results of Connes and Woods [CW] on the classification of such systems, particularly those that are approximately transitive, and by the ordered K$_0$-theoretic classification of C*-algebras initiated by Elliott [E], Giordano and I [GH] constructed a complete invariant for the measure theoretic classification of ergodic $T$, called a dimension space. [This formed part of a more general construction, specialized to $\Z$-actions (that is, single automorphisms).]

We outline the construction, described in detail in [GH]. Let $A = \R[x,x^{-1}]$ be the usual Laurent polynomial ring. We use inner product notation to describe the coefficients of a polynomial $p$, that is, $ p = \sum (p,x^k) x^k$. The ring $A$ is equipped with the obvious partial ordering, $A^+ = \Set{p \in A}{(p,x^k) \geq 0 \text{ for all $k$}}$, making it into a partially ordered ring. We impose $l^1$-norm on $A$, $\| p \| = \sum |(p,x^k)|$, and of course the completion is $l^1 (\Z)$, with convolution extending the multiplication. The evaluation (or augmentation) map, from $A$ to $\R$, given $p \mapsto p(1)$ is of course positive, and $\| p \| \geq |p(1)|$ with equality when $p \in A^+$.

Now form the space of columns of size $n$, denoted $A^n$; equipped with the obvious (coordinatewise) positive cone, $(A^n)^+ = \Set{(p_1, p_2,\dots, p_n)^T \in A^n}{p_i \in A^+}$, $A^n$ becomes a partially ordered $A$-module, and multiplication by $x$ (described by $(p_j)^T \mapsto (xp_j)^T$) is a positive invertible $A$-module transformation. In addition, the augmentation map extends to the obvious map $A^n \to \R^n$, which is again positive.

Now let $n(1), n(2), \dots,$ be an infinite sequence of positive integers, and let $M_k$ be $n(k+1) \times n(k)$ matrices with entries from $A^+$. We can also evaluate each matrix entrywise using the augmentation map; the resulting real matrix $M_i(1)$ (obtained by evaluating each entry at $x=1$) is of course nonnegative, and thus induces an order-preserving map $\R^{n(k)} \to \R^{n(k+1)}$. This gives rise to a (preliminary) direct limit construction,

$$
\diagram
A^{n(1)}&\rTo^{M_{1}} & A^{n(2)} & \rTo^{M_{2}} &
A^{n(3)} & \rTo^{M_{3}} & \dots & & G = \lim_{\to} \Arrow M_k;A^{n(k)}.A^{n(k+1)}\\
\dTo^{} &&\dTo^{}& & \dTo^{} & &&&\dTo^{\rho} \\
\R^{n(1)}&\rTo^{M_{1}(1)} & \R^{n(2)} & \rTo^{M_{2}(1)} &
\R^{n(3)} & \rTo^{M_{3}(1)} & \dots & & G(1) = \lim_{\to} \Arrow M_k(1);\R^{n(k)}.\R^{n(k+1)}\\
\enddiagram
$$
The direct limit of the top row, G, is a partially ordered $A$-module; the direct limit of the bottom row is a partially ordered vector space, and the vertical map(s), augmentation, induce a positive onto map. Now we make another assumption on the matrices, that the column sums of all the $M_i (1)$ are one. [This justifies the alternative name, matrix-valued random walks.]

Now we have to find a linear functional on the direct limit which will (eventually) translate to the measure $\mu$. Elements of the $j$th term in the top row will be indicated by $(v,j)$ where $v \in A^{n(j)}$; the image of the latter in the direct limit will be denoted $[v,j]$ (so that $[M_j v,j+1] = [v,j]$). Because we will be working almost exclusively with invariant measures, we are interested here only in $x$-invariant linear functionals, that is, a positive linear functional $\tau$ on $G$ \st $\tau [xv,j] = \tau[v,j]$. By [GH, Lemma 2.3], these must factor through the vertical map $G \to G(1)$.

The direct limit $G(1)$ is itself a dimension group (over the reals, rather than the usual integers), and under weak conditions (e.g., all the $M_i(1)$ having no zero entries) will itself admit a positive linear functional (known as a trace). Let $\Arrow \phi;G(1).\R$ be a trace on $G(1)$; then $\tau:=\phi \comp \rho$ is a positive invariant linear functional on $G$.

Now we complete $G$ \wrt the L${}^1$ norm induced by $\tau$; this construction can be found in [GoH]. On elements of any partially ordered abelian group $J$ with a trace $\tau$, define the (possibly pseudo)-norm $\| j \| = \inf \Set{\tau(j_1) + \tau(j_2)}{j_1, j_2 \in J^+ \ \& \ j = j_1 - j_2} $. In the case of $G$, the completion is of the form $H:= L^1(Y,\mu)$ and the action of multiplication by $x$ is an invertible positive isometry thereon. (We can thus regard the completion, $H$, as an ordered $l^1(\Z)$-module. Moreover, $\mu$ is determined by the trace $\tau$: $\tau$ extends to the completion, and $\int_Y f \,d\mu = \tau(f)$.

The action of $x$ on $Y$ (induced by the action of $x$ on the $H$) is ergodic if and only if $\rho$, the trace on the vector space $G(1)$, is extremal as a trace thereon [GH, Proposition 3.8]. Although the $Y$ is impossibly complicated to work with directly, we can easily decide when the action will be ergodic.

Not only does this construction yield an ergodic $\Z$-action, but all ergodic $\Z$-actions with invariant measure can be realized by this method of construction (this is part of a more general result of Elliott \& Giordano, that every ergodic amenable action of a discrete group can be constructed using $\R G$ in place of $A = \R[x,x^{-1}] = \R[\Z]$, with suitable Radon-Nikodym derivatives in the non-invariant case).

Although the space $(Y,\mu)$ is normally far too complicated to work with, we can decide isomorphism of two such systems relatively easily. Here isomorphism means that the corresponding $\Z$-actions are measure theoretically conjugate. The ordered $A$-module, together with its associated trace (measure) is known as the {\it dimension space\/} associated to $(M_j)$. If $H$ and $H'$ with all their structures (that is, with the trace, with the ordered $A$-module structure) are isomorphic, then the actions are isomorphic---this is almost tautological. However, when both are ergodic, there is a criterion which allows us to sometimes determine isomorphism.

Suppose $H$ is constructed from the sequence $(M_j)$ and $H'$ is constructed from the sequence $(M_k')$. Each comes with its own L$^1$-norm induced by the ergodic invariant trace (measure). There is a notion of approximate intertwining which is equivalent to isomorphism of the dimension spaces, as follows.

First, we may form telescopings of $(M_j)$ and $(M_k')$, that is, two strictly increasing sequence $1 \leq m(1) < m(2) < m(3) <\dots$ and $1 \leq m'(1) < m'(2) < m'(3) <\dots$ from which we define new sequences $(M^{(j)} = M_{m(j+1)-1}M_{m(j+1)-2}\cdot \dots \cdot M_{m(j+1)}M_{m(j)})$ and $(M^{(k)}{}' = M_{m'(k+1)-1}M_{m(k+1)-2}\cdot \dots \cdot M_{m'(k+1)}M_{m'(k)})$. Now we suppose we have sequences of rectangular matrices $R_k$ and $S_k$ with entries from $A^+$ \st $S_kR_k$ and $M^{(k)}$ have the same dimensions, and $R_{k+1}S_k$ and ${M^{(k)}}{}'$ have the same dimensions. We do not insist on the equalities, $S_kR_k = M^{(k)}$ and $R_{k+1}S_k=M^{(k)}{}'$, but instead that they are close, in the following sense.

We can consider the difference $M^{(k)} - S_kR_k$ as a map $A^{n(k)} \to A^{n(k+1)}$. Each of its columns corresponds to an element of $G$, and we take its norm in $G$ (the L$^1$-norm induced by the trace); then take the maximum over all the columns (this corresponds to the $1-1$ norm, as a map between L$^1$-spaces. Denote the resulting norm via $\tripnorm {M^{(k)} - S_kR_k}xxx$. We similarly define $\tripnorm {R_{k+1}S_k - M^{(k)}{}'}xxx$ (the notation does not reflect the dependance on the choice of dimension space; this would be too cumbersome). Then we say the resulting diagram
$$
\diagram
A^{m(1)}&\rTo^{M^{(1)}} & A^{m(2)} & \rTo^{M^{(2)}} &
A^{m(3)} & \rTo^{M^{(3)}} & \dots & & \\
\dTo^{R_1} &S_1{{\nearrow}} &\dTo^{R_2}&{S_2} {\nearrow}& \dTo^{R_3} & {S_3}{\nearrow}&&& \\
A^{m'(1)}&\rTo^{M^{(1)}{}'} & A^{m(2)} & \rTo^{M^{(2)}{}'} &
A^{m'(3)} & \rTo^{M^{(3)}{}'} & \dots & & \\
\enddiagram
$$
is {\it approximately intertwining\/} if
$$
\sum_k \tripnorm {M^{(k)} - S_kR_k}xxx < \infty \qquad \text{and} \qquad \sum_k \tripnorm { {R_{k+1}S_k - M^{(k)}{}'}}xxx < \infty.
$$
Finally, the two dimension spaces (respectively, the ergodic transformations) are isomorphic if and only if there exists an approximately intertwining diagram between them [GH, Theorem 3.1]. We can also refine this to ensure that the the $R$s and $S$s are norm one. Despite the complications, computations are frequently quite simple, especially in dealing with circulant and related matrices.

Order preserving $A$-module maps between dimension spaces are implememted by the one-sided version of the preceding (no $S$s involved).

The dimension space (equivalently, the transformation, or the sequence of matrices $(M_j)$) is said to be AT$(n)$ if it can be represented as a sequence with matrix sizes at most $n$ (this is not the actual definition, but is equivalent to it, [GH, pp\,32--33]), and it is AT (or {\it approximately transitive\/}) if it is AT(1) (again, this is not the actual definition, but is equivalent to it).

In this paper, we typically deal with dimension spaces arising from sequences $(M_j)$, where the $M_j$ which all have left eigenvector $\pmb 1_n = (1,1,\dots,1)$ (with constant $n$), and we use the resulting $\rho$ to obtain the invariant linear functional $\tau$. We first have to investigate when the this is ergodic, which we do in section 1. We also show that for this class of examples, the matrix of squares, $(M_j^2)$ (and higher, varying powers) is AT. (There are very difficult examples of this type arising from $n=2$ for which $(M_j)$ is AT (but is necessarily AT(2).)\vskip 4pt

\noindent {\it Statement of results.}
Section 1 provides the preliminaries on sequences of hemicirculant matrices (a generalization of circulant matrices, but compatible with products of abelian groups), and a suprising result, that if they are squares (or higher powers), then the resulting sequence, if ergodic, yields an AT action. Section\, 2 deals with a strong property that guarantees isomorphism of the dynamical system with the sequence of traces, called hollowness. We show (as a special case), that if $M_j = \(\smallmatrix 1 & x^{g(n)}\\ x^{g(n)}& 1 \\\endsmallmatrix\)$, then the dynamical system to which it corresponds is hollow if $g$ satisfies an even-term recurrence relation, and is not hollow (but still AT) if it satisfies an odd-term recurrence relation.

Section 3 deals with tensor products. We can construct tensor products of dimension modules in a natural way, thereby obtaining a construction of new dynamical systems (which always preserves ergodicity), it is not clear how to obtain these dynamically (without reference to dimension modules). In any event, under very modest conditions, $(M_j^2)$ and $(M_j \otimes M_j)$ yield conjugate systems, and this is used in a later example.

Section 4 recalls a class of numerical invariants from [H2], and uses it to show that if $(M_j^2)$ is not hollow, then $M_j$ is not isomorphic to its sequence of traces (the latter is automatically AT, the former generically). Combined with the previous section, we obtain a sequence of hemicirculant matrices $(M_j)$ \st $(M_j \otimes M_j^T)$ and $(M_j M_j^T)$ do not yield conjugate systems, in contrast to the main result of section\, 3.

The brief section 5 explains how sequences of hemicirculant matrices arise naturally from the dual action of product type actions of finite abelian groups on some W*-algebras.

Section 6 deals with a general problem; given an ergodic dynamical system $(X,T,\mu)$ with corresponding dimension module, how do we construct the dimension module for a power of $T$, that is, $(X,T^n,\mu)$ (assuming $T^n$ is ergodic). This is very closely related to systems of circulant matrices (suitably modified). Among other things, we show that if $T$ is AT and $T^k$ is ergodic, then $T^k \otimes T^k$ is AT, a somewhat mystifying result (since $T^k$ need not be AT).

\SecT 1 Circulant matrices and their relatives

Let $H$ be a finite abelian group of order $n$, form the group algebra $V = \R H$, treating it as a vector space with basis $\brcs{e_g}_{g \in H}$, on which $\R H$ (a different copy) acts as a commutative algebra of endomorphisms, via
$m_g (e_h) = e_{gh}$. Then each $m_g$ is represented (\wrt the basis $\brcs{e_g}$) as a size $n$ permutation matrix (arising from the regular representation of $H$), and we identify $\sum_{g} q_g m_g$ with the matrix that represents it (where $q_g \in A = \R[x,x^{-1}]$). The matrix representations of $m_g$ and $m_h$ have disjoint supports (that is, for each coordinate, at most one of them has a $1$ in that position), and $\sum_H m_g$ is represented by the matrix all of whose entries are $1$.

For $ M = \sum_H q_g m_g$ a matrix with entries from $A$, then all of its entries lie in $A^+$ if and only if all belong to $q_g $ to $A^+$ (this comes from disjointness of the supporting entries); the collection of such will be denoted $(A H)^+$ (viewed as a subset of the $n \times n $ matrix algebra over $A$).

We can write down explicitly the common eigenvectors for the elements of $A H$; for each $\alpha $ in $\hat H$ (the dual group of $H$), define a vector $v_{\alpha}$ in $V \otimes \C = \C H$ via
$$
v_{\alpha} = \frac 1{|H|} \sum_{g \in H} \alpha(g) e_g
$$
(the normalization is to ensure that $\| v_{\alpha}\| = 1$ in the appropriate $l^1$-like norm). We check immediately that $m_g v_{\alpha} = \alpha(g^{-1}) v_{\alpha}$, so that each $v_{\alpha}$ is an eigenvector, they are distinct (and can be separated by the $m_g$), and form an orthonormal set, hence constitute a basis for $\C H$. For any $M = \sum_H q_g m_g$ in $A H$, its $\alpha$th eigenvalue, $\lambda_{\alpha}(M)$, is in $\C[x,x^{-1}]$, and is given by the eigenvalue corresponding to $v_{\alpha}$, that is,
$$
\lambda_{\alpha} (M) = \sum_{g \in H} q_g \alpha(g^{-1}).
$$
For each $\alpha$, the assignment $\Arrow \lambda_{\alpha}; A H . \C[x,x^{-1}]$
is an $A$-algebra homomorphism (that is, addititive, multiplicative, and compatible with multiplication by elements of $A$). A special choice occurs if $\alpha$ is the trivial character, denoted $\chi_0$, in which case we denote the corresponding eigenvector $v_0 = \frac 1{|H|} \sum e_g$, and $\lambda_0 (M) = \sum q_g$.

This defines right eigenvectors; left eigenvectors, $w_{\alpha}$, are defined analogously, but with complex conjugation. The corresponding left eigenvector is given by $w_{\alpha} = \sum E_g \alpha(g^{-1})$ where we are using $E_g$ to denote basis elements of $AH$ viewed as a {\it right\/} $AH$ module. In particular, $w_{\alpha} = |H| \overline{v_{\alpha}}^T$ (note the complex conjugate); by construction, $w_{\alpha} v_{\alpha} = 1$, and as an operator $\Arrow w_{\alpha} ; AH.A$ (via $w_{\alpha} (e_g) = \alpha(g^{-1})$), in the operator $l^1$-norm, $\| w_{\alpha} \| = 1$.

\Lem Examples.

\noindent (a) {\it Circulant matrices\/}.
If $H = \Z_n$, the cyclic group of size $n$, then $AH$ consists precisely of the circulant matrices of size $n$ (with entries from $A$), where the generator $g = [1]$ maps to the cyclic permutation matrix whose first row is $(0 \ 1 \ 0 \dots \ 0)$.

\noindent (b) $H = \Z_2 \times \Z_2$. In this case, $ n = 4$, and $AH$ consists of all matrices of the form
$$
\(\matrix
a_0 & a_1 & a_2 & a_3 \\
a_1 & a_0 & a_3 & a_2 \\
a_2 & a_3 & a_0 & a_1 \\
a_3 & a_2 & a_1 & a_0\\
\endmatrix\) = a_0 \I + a_1 \(\matrix
0 & 1 & 0 & 0 \\
1 & 0 & 0 & 0 \\
0 & 0 & 0 & 1 \\
0 & 0 & 1 & 0\\
\endmatrix\) + a_2 \(\matrix
0 & 0 & 1 & 0 \\
0 & 0 & 0 & 1 \\
1& 0 & 0 & 0 \\
0 & 1 & 0 & 0\\
\endmatrix\) + a_3 \(\matrix
0 & 0 & 0 & 1 \\
0 & 0 & 1 & 0 \\
0& 1 & 0 & 0 \\
1 & 0 & 0 & 0\\
\endmatrix\)
$$
where $a_i $ belong to $A$.
\qed

We call a matrix {\it hemicirculant \wrt $H$} (or $H$-{\it hemicirculant\/}) if it is of the form $M = \sum_H q_g m_g$. Suppose that $(M_j)$ is a sequence of hemicirculant nonnegative matrices (\wrt the same, fixed, $H$), that is,
in $ (AH)^+$, \st $\lambda_0 (M_j) (1) = 1$, so that the eigenvalue at $v_0$, when evaluated at $x = 1$, is $1$, i.e., $\sum q_g(1) = 1$. We wish to study the dimension space arising out of the direct limit $\lim \Arrow M_j ; A^n . A^n$, where we take the invariant functional arising from the left eigenvector $(1 \ 1 \ \dots \ 1)$ followed by evaluation at $x = 1$.

We must first check when this yields an ergodic dimension space. This boils down to weak ergodicity (in the context of matrices with real entries) of hemicirculant matrices, which should be in the literature, but I was not able to find any reasonable references (except of course when $n $ is prime).

Let $C_j$ be hemicirculant matrices with nonnegative {\it real\/} entries, and whose column sums are all $1$. Weak ergodicity of the sequence $(C_j)$ is equivalent to the limiting dimension group $\lim \Arrow C_j ; \R^n .\R^n$ having unique trace, and also to projective and actual convergence of the products to a rank one matrix (this is the usual definition, although in general, order matters; here, the matrices mutually commute), which in this case must be $v_0 w_0$.

Write $C_j = \sum_{g \in H} c_{gj} m_g$ where $0 \leq c_{gj}$ and $\sum_{g \in H} c_{gj} = 1$ for all $j$. It is easy to check that if $n$ is prime (so in particular, the $C_j$ are circulant matrices), then weak ergodicity is equivalent to $\sum_j \(1 - \max_{g}\brcs{c_{gj}}\) = \infty$. For other choices of $H$ (even cyclic ones), this criterion is merely necessary, but not sufficient (the mass may accumulate on a proper subgroup of $H$). While the criterion becomes messier and messier depending on how far from cyclic $H$ is, it is still computable, and more importantly, it is used in our subsequent AT results.

Fix $H$, and let $\alpha$ and $\beta$ be two unequal elements of $\hat H$. Define $S_{H,\alpha, \beta}$ (or $S_{\alpha, \beta}$ if there is no ambiguity about $H$) via
$$
S_{H,\alpha,\beta} = \Set{(g,h) \in H \times H }{\alpha(g^{-1}) \beta(h^{-1}) \neq \alpha (h^{-1})\beta(g^{-1})}.
$$
Since $\alpha (g^{-1})$ is just the complex conjugate of $\alpha (g)$, we could have just written the defining condition as $\alpha(gh^{-1}) \neq \beta (gh^{-1})$, or even $\alpha \beta^{-1} (gh^{-1}) \neq 1$; in particular, whether $(g,h)$ belongs to $S_{\alpha,\beta}$ depends only on $gh^{-1} $ and $\alpha\beta^{-1}$.

The statement in the following criterion involves a centring (otherwise it becomes even more awkward); the condition that $c_{0j}$ be maximal can be arranged by multiplying each term by a suitable element of $H$; this does not affect weak ergodicity. The criteria in (iii) boil down to one computation for each prime divisor of $n = |H|$ if $H$ is cyclic, but hordes of them in the noncyclic case (e.g., if $H = (\Z_p)^k$ with prime $p$, then the number of maximal subgroups is $p^k - 1 = n-1$). It may be possible to restrict the maximal subgroups to a reasonable number.

\comment
\Lem Lemma. (Weak ergodicity criterion) Suppose for each $j$ that $C_j = \sum_{g \in H} c_{gj} m_g$ is a real hemicirculant matrix \wrt\ $H$, where $0 \leq c_{gj}$ and $\sum_{g \in H} c_{gj} = 1$. Assume in addition that $c_{0j} = \max_{g \in H} \brcs{c_{gj}}$ for all $j$. The following are equivalent.
{\par} \item{(i)} The limit dimension group $\lim \Arrow C_j; \R^n . \R^n$ is a simple dimension group with unique trace.
{\par} \item{(ii)} For each maximal proper subgroup $K$ of $H$, $\sum_j \sum_{g\not\in K} c_{gj} = \infty$
{\par} \item{(iii)}
{\par} \item{(iv)} For all $\alpha \neq \beta$ in $\hat H$, $\sum_j \sum_{(g,h) \in S_{\alpha,\beta}} c_{g j}c_{h j} = \infty$.

\endcomment

\Lem Lemma \oneone. (Weak ergodicity criterion) Suppose for each $j$ that $C_j = \sum_{g \in H} c_{gj} m_g$ is a real hemicirculant matrix \wrt\ $H$, where $0 \leq c_{gj}$ and $\sum_{g \in H} c_{gj} = 1$. Assume in addition that $c_{0j} = \max_{g \in H} \brcs{c_{gj}}$ for all $j$. The following are equivalent.
{\par} \item{(0)} The limit dimension group $\lim \Arrow C_j; \R^n . \R^n$ is a simple dimension group with unique trace;
{\par} \item{(i)} for every $\alpha \neq \chi_0$, for all $j_0$, $\prod_{j \geq j_0} \lambda_{\alpha} (C_j) = 0$;
{\par} \item{(ii)} for every $\alpha \neq \chi_0$, for all $j_0$, $\sum_j \sum_{\alpha(g) \neq 1} c_{gj} = \infty$;
{\par} \item{(iii)} for each maximal proper subgroup $K$ of $H$, $\sum_j \sum_{g\not\in K} c_{gj} = \infty$;
{\par} \item{(iv)} for all $\alpha \neq \beta$ in $\hat H$, $\sum_j \sum_{(g,h) \in S_{\alpha,\beta}} c_{g j}c_{h j} = \infty$.

\Pf The equivalence of (0)--(ii) is clear, and the equivalence of (ii) with (iii) stems from every proper subgroup being contained in a maximal proper subgroup. Let $n = |H|.$

We prove (ii) implies (iv). Form the element $\gamma = \alpha \beta^{-1}$ of
$\hat H$; this is not the trivial character. For each $j$, we may find
$g_j$ not in the kernel of $\gamma$ \st $c_{g_j j} \geq \sum_{\gamma(g)
\neq 1} c_{gj} /n$; this forces $\sum_j c_{g_j j} = \infty$, and since $c_{0,j}
\geq 1/n$ as a consequence of the hypotheses, we have $\sum_j c_{g_j j}
c_{0,j}= \infty$.

We note that $(g,h) \in S_{\alpha,\beta}$ if $gh^{-1}$ does not belong
to the kernel of $\gamma$. Hence $(g_j, 1) \in S_{\alpha, \beta}$, and
thus $\sum_j \sum_{g \in S_{\alpha,\beta}} c_{gj}c_{h,j} \geq \sum_j
c_{g_jj}c_{0,j}$, and the latter diverges.

Now (iv) implies (ii). Let $\gamma$ be a nontrivial element of $\hat H$. With
$\alpha = \gamma$ and $\beta = \chi_0$ (the trivial character), find
$(g_j, h_j) \in S_{\alpha,\beta}$ \st $c_{g_j j} c_{h_j j} \geq
\sum_{S_{\alpha,\beta}} c_{gj} c_{hj}/n^2$ (this is possible since the sum
is over fewer than $n^2 $ elements). Since $\gamma(gh^{-1}) \neq 1$, at
least one of $c_{g_j}$ or $c_{h_j}$ is not in the kernel of $\gamma$. We are
done, since $ \sum_j \sum_{\gamma (g) \neq 1} c_{gj} \geq c_{g_j}$ or
$c_{h_j}$ respectively. \qed

For elements $\alpha$ and $\beta$ of $\hat H$, and elements $g$ and $h$, set $z =\alpha(g^{-1}) \beta(h^{-1})$ and $y =\alpha (h^{-1})\beta(g^{-1})$. Let $N$ denote the exponent of $H$ (smallest positive integer exceeding $1$ \st $k^N = 1$ for all $k$ in $H$), and set $\xi = e^{2\pi i/N}$. Then the values of all irreducible characters on elements of $H$ lie in $\brcs{\xi^j}$. Obviously, $|z + y| \leq 2$, with equality occurring only if the arguments of $z$ and of $y$ are equal---in this case implying $z = y$ (since both are just powers of $\xi$). On the other hand, if $|z + y | < 2$, then $ |z+y| < 2 \cos \pi/N$ (not $2\cos 2\pi/N$). In particular, either $z = y$ or $|z + y| \leq 2 \cos \pi/N$, and so if $(g,h)$ belongs to $S_{H,\alpha,\beta}$, then $|\alpha(g^{-1}) \beta(h^{-1}) + \alpha (h^{-1})\beta(g^{-1})| \leq 2 \cos \pi/N$.

When we write $(M_j)$ is an ergodic sequence of hemicirculant matrices (\wrt some finite group $H$), we mean that $M_j \in (AH)^+$, that the column sums, on evaluation at $x = 1$, are all $1$, and the sequence $(C_j = M_j(1))$ is weakly ergodic; in particular, the trace obtained from the constant eigenvector yields an ergodic measure. The eigenvalues of $M_j$ lie in $\C[x,x^{-1}] = A \otimes \C$; we impose the obvious $l^1$-norm on the latter (sum of the absolute values of the coefficients).

\comment
new stuff

Let $(A_i)$ be a sequence of $p\times p$ circulant matrices with entries in $S^+$ \st the system $(A_i)$ is ergodic. Let $v_i = (\zeta^{ij})$ (indexed modulo $p$ in both coordinates), and set $\Lambda_2 (A) $ to be the maximum $l^1$-norm of all the eigenvalues of $A$ other than the row sum; in particular, $\Lambda_2 (A) \leq 1$ with strict inequality except in some degenerate instances. In other words, if $A = \sum_j p_j(X)B^j$, then $\Lambda_2 (A) = \max_{k\neq 0} \brcs{\| \sum_{j=1}^{p-1}
p_j (X) \zeta^{kj} \|}$.

\Lem Proposition. If $\sum_i \( 1 - \Lambda_2 (A)\) = \infty$, then $(A_i)$ is AT and given by the sequence of row sums.

\Pf First work with the complexified version of $S$. We have a natural map $\psi_n;S^p .S$ given by left multiplication by $(1 \ 1 \dots 1) = v_0^T$. We show that for all $v_j$ with $j \neq 0$, given $\epsilon$ and $m$, there exists $M \equiv M(\epsilon,m)$ \st $\|A_{n+M} A_{n+M-1}\dots A_{n+1} v_j \| < \epsilon$, and this will still hold if $M$ is increased. The expression within the norm signs is simply $\(\prod_{k=1}^M \lambda_j(A_{n+k})\) v_j $, where $\lambda_j$ is the eigenvalue corresponding to $v_j$.
We have
$$\eqalign{
\left\|\prod_{k=1}^M \lambda_j(A_{n+k})\right\| & \leq \prod_{k=1}^M \left\| \lambda_j(A_{n+k})\right\| \cr
& \leq \prod_{k=1}^M \Lambda_2(A_{n+k}) \cr
}$$
We can write $\Lambda_2 (A_{n+k}) = 1 - (1-\Lambda_2 (A_{n+k}))$, so the hypothesis implies that the product goes to zero.

Now $\brcs{v_j}$ is an $S$-basis for $S^p$. Hence given $\epsilon$ and $n$, there exists $N \equiv N(n, \epsilon)$ \st for any $w$ written uniquely as $w = \sum f_i(w) v_{i}$ in $S^p$, $\|\prod_{k=1}^N A_{n+k} w - \prod \lambda_0 (A_{n+ k}) f_0(w_0) v_0 \| < \epsilon \| w\|$.

This permits us to construct an approximately commuting isomorphism between the two systems $(A_i) $ and $(\tr A_i)$; first, the vertical map downward from $S^p$ to $S$ given by $\psi_n$. Now let $\epsilon (j)$ be a summable sequence of $\epsilon$s, and let $N(1) < N(2) < \dots$ be defined so that having defined the integer $N(j)$, we have $\|\prod_{k=N(j)+1}^{N(j+1)} A_{k} w - \prod \lambda_0 (A_{k}) f_0(w_0) v_0 \| < \epsilon \| w\|$. Then the maps $\Arrow \phi_j; S . S^{p}$ defined from the $N(j)$ copy of $S$ to the $N(j+1)$ copy of $S^p$ approximately commute and it is easy to see that this yields an isomorphism. Moreover, both families of maps restrict to the real versions of $S$, $S^p$, and they are positive, and it easily follows that they induce mutually inverse maps here as well, yielding an isomorphism at this level.
\qed

Let $(A_i)$ be a sequence of $p\times p$ circulant matrices with entries
in $S^+$ \st the system $(A_i)$ is ergodic. Let $v_i = (\zeta^{ij})$
(indexed modulo $p$ in both coordinates), and set $\Lambda_2 (A) $ to be
the maximum $l^1$-norm of all the eigenvalues of $A$ other than the row
sum; in particular, $\Lambda_2 (A) \leq 1$ with strict inequality except
in some degenerate instances. In other words, if $A = \sum_j p_j(X)B^j$,
then $\Lambda_2 (A) = \max_{k\neq 0} \brcs{\| \sum_{j=1}^{p-1}
p_j (X) \zeta^{kj} \|}$.

\endcomment
\comment
\Lem Prop. Let $\brcs{A_j}$ be a sequence of circulant matrices, each
written as
$A_j = \sum_{i} q_{ij}(x) B^i$ with $q_{ij} \in S^+$ and $\sum_i q_{ij}(1)
=1$ for all $j$, and suppose that $\(\Arrow A_i; S^p . S^p\)$ is an
ergodic system. Suppose in addition that each $A_j$ is symmetric.
Sufficient for the system to be AT is that $\sum_{j} (1- \sum_{(i,p) > 1}
q_{ij}(1)) = \infty$.

\Pf By symmetry, we may write $A_{j} = p_{0j} + \sum_{i=1}^{\flo{p/2}}
s_{ij}(B^{i} + B^{-i})$, where $2s_{ij} = p_{ij} + p_{-i,j}$. With $k \neq
0$, $A_j v_k = v_k (p_{0j} + \sum s_{ij} (\zeta^{ki} + \zeta^{-ki})$, so
the eigenvalue, denoted $e_{kj}$ is $p_{0j} + 2\sum_{i=1}^{\flo{p/2}}
s_{ij}\cos 2ki\pi/p $. Hence the $l^1$ norm of the eigenvalue is bounded
above by $p_{0j}(1) + 2\cos \pi/p\sum_{i \in J_k} s_{ij}(1) + 2 \sum_{i
\in J_k^c} s_{ij}(1)$, where $J_k = \Set{1 \leq i \leq \flo{p/2}}{p \text{
does not divide }ik}$, and its complement, $J_k^c$, is taken \wrt the set
$\brcs{1,2,\dots, \flo{p/2}}$. Since $1 = \sum_i q_{ij} (1)$, we deduce
that
$$
\left\| e_{kj} \right \| \leq 1 - (2 - 2\cos \pi/p)\sum_{i \in J_k}
s_{ij}(1)
= 1 - (2\sin^2 \pi/2p) \sum_{i \in J_k} (2s_{ij}).
$$
If $p$ is prime, then this is simply $1- 2\sin^2 \pi/2p (1 - p_{0j} (1))$,
hence if $\sum_j (1- p_{0j} (1)) = \infty$, then $\left\|\prod_{j \geq N}
e_{kj} \right\| \to 0$ for all $N$. This would force the element $(v_k,
N)$ (in the $N$th copy of $S^p$ to go to zero in the system. Hence the
obvious maps are going to be isomorphisms.

[basically an observation that the eigenvalues other than the row sum get
smaller in norm, converging to zero in the products, so contribute nothing
to the kernel of the vertical map. More general than it seems since
telescoping can improve the thing to be summed! e.g., $p= 6$, each $A_j$
has $q_{ij} = 0$ for $i= 1,4,5$; telescope in
pairs, and if there is sufficient mass at $q_{2j}$ and $q_{3j}$, then get
mass at $i=5$ in the telescoped product. ]

Similar but easier result with $AA^T$.

Lemma. Suppose that $\brcs{c_i}_{i=1}^{k}$ is a bunch of nonnegative real
numbers adding to one, and let $\alpha = \max{c_i} \geq 1/k$. Let $l =
\flo{1/\alpha}$. Then
$$
\min \sum_{i < j} c_i c_j =
\cases
\alpha (1- \alpha) & \text{if $\alpha \geq 1/2$} \\
l\alpha (1 - (l+1)\alpha/2) & \text{if $\alpha \leq 1/2$}.
\endcases
$$

\Pf Consider the constrained problem, to minimize $\sum_{i < j} X_i X_j$
subject to $X_i \geq 0$, $\sum c_i = 1$, and $X_i \leq \alpha$ for all
$i$. If we look at the interior of the feasible domain (that is, $\alpha>
c_i > 0$ and $\sum X_i = 1$, on setting $X_k = 1 - \sum_{i=1}^{k-1} X_i$,
we see that there are no critical points of interest for this problem, in
fact, at any critical point, we must have at least two zero coordinates.
Remove all the zeros, and we are left with at least one of the values
hitting $\alpha$ (at an optimal point). Say there are $t$ of the $c$s
which hit $\alpha$; necessarily $t \leq l+1$ (and if $\alpha$ is
irrational, $t \leq l$). By relabelling, we may assume $c_k = \alpha =
c_{k-1} = c_{k-t+1}$. This reduces the problem to minimizing
$$
\sum_{i < j < k-t+1} X_i X_j + (1-t\alpha)t\alpha + (t^2 -t) \alpha^2/2,
\text{subject to $\alpha> X_i > 0$, $\sum_{i < k-t+1} X_i = 1 - t
\alpha$.}
$$
Again, the only time a critical point can occur is when there is at most
one among the $X_i$ that is not zero. Necessarily, either the nonzero
entries are simply $(\alpha,\alpha, \dots, \alpha, 0, 0, \dots, 0)$ ($l$
$\alpha$s) or there are $k$ of the $\alpha$s, and the remaining term is
$1-l \alpha$. If $\alpha > 1/2$, the minimal value occurs at
$(\alpha,1-\alpha,0,0,\dots,0)$ and is $\alpha(1-\alpha)$. If $\alpha \leq
1/2$ (so $l \geq 2$), the critical point is either $(1/l,\dots,
1/l,0,\dots,0)$ (if $\alpha = 1/l$) or $(\alpha,\alpha, \dots, \alpha,
1-l\alpha, 0, \dots, 0)$ with $l$ of the $\alpha$s. In the former case,
the value is $(l^2 -l)/2l^2 = 1/2 - 1/2l = (1 - \alpha)/2$; in the latter
case, it is $(1-l\alpha) l \alpha + (l^2 -l) \alpha^2/2 = l\alpha (1 -
(l+1)\alpha/2)$.
\qed

Note in particular, that when $\alpha \leq 1/2$, the minimal value is at
least $1/4$, converging to $1/2$ as $\alpha \to 0$ (this requires $k \to
\infty$), and it is monotone decreasing in $\alpha$ on $(0,1)$. We can
rephrase this as $\sum c_i c_j \leq h(\alpha) \cdot (\sum c_i)^2$ where
$h(\alpha)$ is given in the display above.

For $AA^T$. Write $A = \sum q_i B^i$, so the $k$th eigenvalue of $A$,
$e_k$, is $\sum_{i} q_i \zeta^{ik}$, and that of $AA^T$ is $e_k
\overline{e_k} = \sum_{i} q_i \zeta^{ik} \sum_{i} q_i \zeta^{-ik}$,
which expands as
$$\eqalign{
\sum_i q_i^2 + \sum_{i < j} q_i q_j (\zeta^{(i-j)k} + \zeta^{(j-i)k}) & =
\sum_i q_i^2 + 2 \sum_{i < j} q_i q_j (\cos 2(j-i)k\pi/p); \text{ hence
}\cr
\| e_k \overline{e_k} \| & \leq \sum_i q_i^2 (1) + 2 \cos \pi/p \sum_{J_k}
q_iq_j(1) + \sum_{J_k^c} q_iq_j(1) \cr
& = 1 - (2-2\cos \pi/p) \sum_{J_k} q_iq_j(1)
}$$
where $J_k = \Set{(i ,j)}{p \text{does not divide $2(j-i)k$}}$.

\endcomment


\comment

\def\gcd{\text{gcd\,}}

As usual, $B$ is the cyclic permutation matrix of size $p$. The following
ergodicity criterion for circulant random walks on a set of size $p$ is
probably well known, but as usual, i couldn't find it. When $p$ is prime,
the usually-stated criterion is $\sum_j (1-\max_i \brcs{c_{ij}}) =
\infty$. For nonprime $p$, the criterion is more complicated, and simplify
it a bit, we first multiply each of the matrices by a suitable power of
$B$ (which obviously does not affect ergodicity) so that the coefficient
of $B^0$ is maximal.

For $p$ a positive integer, and $a$, $b$ integers, define $S_{a,b}:=
\Set{(s,t)\in \Z_p \times \Z_p}{(a-b)(s-t) \not\equiv 0 \mod p}$.
Obviously, $S_{a,b} = S_{a-b,0}$ and if $u$ is relatively prime to $p$,
then $S_{a,b} = S_{a',b'}$ whenever $a-b = t(a'-b')$. Finally, $\cap
S_{a,b} = \Set{(s,t)}{\gcd (i-j,p) = 1}$.

\Lem Lemma. (Ergodicity criterion) Suppose $C_j = \sum_{i \in \Z_p}
c_{ij}B^j$ with $c_{ij} \geq 0$, and for all $j$, $\sum_i c_{ij} = 1$.
Assume in addition that $c_{0j} = \max_i\brcs{c_{ij}}$ for all $j$. The
following are equivalent.{\par}
\item{(i)} The limit dimension group $\lim \Arrow C_j; \R^p.\R^p$ is
simple with unique trace.
\item{(ii)} For each prime $r$ dividing $p$, $\sum_j \sum_{\gcd(i,p) = 1}
c_{ij} = \infty $.
\item{(iii)} For each prime $r$ dividing $p$, $\sum_j \sum_{(s,t)\in
S_{p/r,0}} c_{sj} c_{tj} = \infty $.
\item{(iv)} For all $(a,b) \in (\Z_p \times \Z_p) \setminus\brcs{(0,0)}$,
$\sum_j \sum_{(s,t)\in S_{a,b}} c_{sj} c_{tj} = \infty $.

\Pf (i) iff (ii) is straightforward, (iv) implies (iii) is trivial, and
(i) implies (iv) is routine. It remains to show (ii) iff (iii). Since
$\sum_i c_{ij} = 1$, it follows from the hypotheses that $c_{0j} \geq
1/p$.

\noindent (ii) implies (iii). Fix the prime divisor $r$ of $p$. For each
$j$, there exists $i_j $ in $\Z_p$, depending on $r$, $p$, and $j$, with
$\gcd (i_j,r) = 1$ \st $c_{i_j,j} \geq \frac 1p \sum_{\gcd (i,r) = 1}
c_{ij}$. Then
$$
\sum_{S_{p/r,0}} c_{sj} c_{tj} \geq c_{i_j,j} c_{0,j} \geq \frac 1{p^2}
\sum_{\gcd (i,r) = 1}.
$$
Thus (ii) implies (iii).

\noindent (iii) implies (ii). Since $|S_{p/r,0}| < p^2$, for each $j$,
there exists $(s_j,t_j)$ in $S_{p/r,0}$ \st
$$
c_{s_j,j}, c_{t_j,j} \geq c_{s_j,j}\cdot c_{t_j,j} \geq \frac 1{p^2}
\sum_{S_{p/r,0}} c_{sj}c_{tj}.
$$
However, at least one of $\brcs{s_j, t_j}$ is relatively prime to $r$;
hence for every $j$,
$\sum_{\gcd(i,r) =1} c_{ij} \geq \frac 1{p^2} \sum_{S_{p/r,0}}
c_{sj}c_{tj}$. Therefore (iii) implies (ii).
\qed

As a consequence of the definitions, $S_{k,k'} = \Set{(l,m) \in \Z_p \times
\Z_p}{\xi^{kl + k'm} \neq \xi^{k'l + km}}$.
\endcomment

\Lem Lemma \onetwo. Suppose that $\(M_j = \sum_{g \in H} q_{g j}m_g\)$ is an ergodic
sequence of $n \times n$ hemicirculant matrices \wrt the finite abelian group $H$. Then for all $\alpha \neq \beta$ in $\hat H$ and
all $j_0 \geq 1$,
$$
\lim_{d\to \infty}\left\|\lambda_{\alpha} \(\prod_{j= j_0}^{j_0
+d}M_j\)\lambda_{\beta} \(\prod_{j= j_0}^{j_0 +d} M_j\)\right\| = 0
$$

\Pf Temporarily drop the subscript $j$. Then
$$\eqalign{
\lambda_{\alpha} (M) \cdot \lambda_{\beta}(M) & = \sum_{g\in H} q_g \alpha(g^{-1})\cdot
\sum_{h\in H} q_h \beta(h^{-1})\cr
& = \frac 12\sum_{(g,h) \in S_{\alpha,\beta}} q_g q_h \cdot (\alpha(g^{-1})\beta(h^{-1}) + \alpha(h^{-1})\beta(g^{-1})) + \sum_{(H\times H)\setminus S_{\alpha,\beta}} \text{\hglue -2.5ex} q_{g}q_h \alpha(g^{-1})\beta(h^{-1});\quad{\text{so}} \cr
\left\| \lambda_{\alpha} (M) \cdot \lambda_{\beta}(M)\right\| & \leq \cos \frac{\pi}{N}
\sum_{(g,h) \in S_{\alpha,\beta}} q_g q_h (1) + \sum_{(H \times H)\setminus
S_{\alpha,\beta}} q_{g}q_h (1) \cr
& = \(\sum_g q_g(1)\)^2 - \(1-\cos \frac{\pi}N\) \sum_{(g,h) \in S_{\alpha,\beta}} q_g q_h
(1)\cr
& = 1 - 2 \sin^2 \frac{\pi}N \sum_{(g,h) \in S_{\alpha,\beta}} q_g q_h
(1). \cr
}$$

Restoring the $j$, by the earlier lemma, with $c_{gj} = q_{gj}(1)$, ergodicity implies that $\sum_j\sum_{S_{\alpha,\beta}} q_{gj} q_{hj} (1) =
\infty$. Thus $\prod\left\| \lambda_{\alpha} (M_j)\lambda_{\beta} (M_j) \right\| \to
0$ (in the strongest possible sense). Since $\lambda_{\alpha}$ is multiplicative
on hemicirculant matrices, the result follows.
\qed

\Lem Theorem \onethr. Suppose that $(M_j)$ is an ergodic sequence of $n \times n$
hemicirculant matrices. Let $\Arrow f;\N.\N$ be a function \st $f(j)\geq 2$ for all $j$. Then $(M_j^{f(j)})$ is AT.
In particular, let $0 = n(1) < n(2) < \dots $ be an infinite sequence of positive integers \st for all $\alpha \neq \beta \in \hat H$, $\sum_i \left\| \lambda_{\alpha}(M^{(i)})\lambda_{\beta}(M^{(i)}) \right\| < \infty$, where $M^{(i)} = \prod_{j=n(i)}^{n(i+1)-1} M_j$. Write $M^{(i)}M^{(i+1)} = \sum_{g\in H} p_{gi}m_g$.{\par}
\item{(a)} Then $(M_j^2)$ is isomorphic to the AT sequence $(\tr (M^{(i)}M^{(i+1)}) =|H| p_{0i})$.
\item{(b)} The automorphism of the dimension space of $(p_{0i})$ induced by $m_g$ (for $g \in H$) is implemented by the map $[f,i] \mapsto [p_{g^{-1}i}f,i+1]$ (for $f \in A$).

\Pf If $M$ is a hemicirculant matrix \st $\left\| \lambda_{\alpha}(M)
\lambda_{\beta}(M)\right\| < \epsilon/p^2$ (as would arise as a product of suitable $M_j$, by the previous lemma), write $M = \sum_{\hat H} \lambda_{\alpha}
(M)v_{\alpha} w_{\alpha} = \sum_g q_g m_g$. Define $V = \sum \lambda_{\alpha}(M) v_{\alpha} $ and $W = \sum \lambda_{\alpha}(M) w_{\alpha} $.We note
$$\eqalign{
\| M^2 - VW \| & =\left\| \lambda_{\alpha} (M)^2\sum v_{\alpha} w_{\alpha} - \(\sum \lambda_{\alpha}(M)v_{\alpha} \)\(\sum
\lambda_{\beta}(M)v_{\beta}\) \right\|\cr
& =
\left\| \sum_{{\alpha}\neq {\beta}} \lambda_{\alpha}(M)
\lambda_{\beta}(M)v_{\alpha} v_{\beta} \right\| \cr
& < \frac{(p^2 - p)\epsilon}{p^2} < \epsilon.
\cr
}$$

From the displayed computation, $\| M^2 - VW\| < \epsilon$. Now we show
that both $V$ and $W$ have entries in $A^+$. Since $W = |H| \overline {V^T}$, it
suffices to show this for $V$.
$$\eqalign{
V &= \sum \lambda_{\alpha}(M) v_{\alpha} \cr
& = \frac 1{|H|} \sum_{\alpha}\sum_g e_g \alpha(g)\lambda_{\alpha} (M) = \frac 1{|H|}\sum_{{\alpha},g,h}
e_g \alpha (g) \alpha(h^{-1})q_h \cr
& = \frac 1{|H|}\sum_{{\alpha},g,h}
e_g q_h \alpha (gh^{-1})\cr
& = \frac 1{|H|}\sum_g e_g\sum_{h} q_h \sum_{\alpha} \alpha (gh^{-1})\cr
&= \sum_g e_g q_{g} \in (AH)^+.\cr
}$$
A similar computation yields $W = \sum q_{g} e_{g^{-1}}$. Since $(M_j)$ is ergodic, for all $j_0$, there exists $d$ \st $\max_{\alpha\neq
{\beta}} \left\|\lambda_{\alpha} (\prod_{j_0}^{j_0+d}M_j) \lambda_{\beta}
(\prod_{j_0}^{j_0+d}M_j)\right\| < \epsilon/p^2$, and it follows that
$(M_j^2)$ is AT.

For higher powers ($f(j)$), the argument is similar.

Now this process can be done for any telescoping $\brcs{n(i)}$ \st $\sum_i \left\| \lambda_{\alpha}(M^{(i)})\lambda_{\beta}(M^{(i)}) \right\| < \infty$ (that is, $M^{(i)} = M$); such a telescoping exists by the previous lemma. Then define $V^{(i)}$ and $W^{(i)}$ as above; the approximate factorization yields an isomorphism between the dimension spaces of $(M_j^2) = ((M^{(i)})^2)$ and $(W^{(i+1)}V^{(i)})$, with maps given by $[z,i]\mapsto [W^{(i)}\cdot z,i]$ ($z \in A^n$) and in the reverse direction via $[f,i] \mapsto [fV^{(i)},i+1]$ ($f \in A$). We expand $W^{(i+1)}V^{(i)} = \frac1{|H|}\sum q_{gi}q_{g^{-1},i+1}$; this is exactly $\frac 1{|H|} \tr M^{(i)}M^{(i+1)} = p_{0i}$.

To determine the effect of the automorphism induced by $g$ in $H$, we begin with the element $[1,i]$ in the dimension space of $(p_{0i})$; under the map
$V^{(i)}$, this is sent to $[V^{(i)},i+1]$ (of the dimension space of $(M_j^2)$); now $g$ acts on this directly by multiplication, yielding the element $[\sum_{h \in H} q_{hi} e_{hg},i+1]$; then $W^{(i+1)}$ sends this to $[\sum_{h \in H} q_{hg^{-1},i}q_{h^{-1},i+1},i+1]$; but this is just $[p_{g^{-1}i},i+1]$. Since the map at the $i$th level is uniquely determined by its effect on $[1,i]$ (and the automorphism exists by the isomorphism), the automorphism is induced as indicated.
\qed

In particular, the diagonal entries of $VW$ are respectively $q_{g^{-1}} q_g$; this is unsurprising, as $\tr M^2 = |H|\sum_g q_g q_{g^{-1}} $ (this entails that $|H|q_g q_{g^{-1}}$ are close in $l^1$ as $g$ varies). The $(g,h)$ entry of $VW$ is $ q_g q_{h^{-1}}$, so $VW$ is not generally itself a hemicirculant matrix and $M$ does not (usually) commute with $VW \sim M^2$; it does not even seem possible to perturb $V$ to $V'$ and $W$ to $W'$ (with error bounded by a fixed multiple of $\sum_{\alpha \neq \beta} \| \lambda_{\alpha}(M) \lambda_{\beta}(M)\|$) so that $V'W'$ commutes with $M$.

\comment
If we attempt the same process for $M_j M_j^T$ (instead of $M_j^2$), not much is gained---since $\lambda_{\alpha}(M) = \lambda_{\alpha^{-1}}(M^T)$, the condition on the distinct dual elements $\alpha$ and $\beta$ includes
$\lambda_{\alpha}(\prod M_j^2) \to 0$ (for $\alpha$ not squaring to the trivial character); this almost implies $(M_j^2)$ is hollow (hollowness, a particularly strong property, will be defined later). It is an easy consequence that if $(M_j^2)$ is hollow, then $(M_j M_j^T)$ is AT, but this is hardly earthshaking (and if $H $ is an elementary $2$-group, is a tautology). On the other hand, if $|H|$ is of odd order, then $(M_j M_j^T)$ is already hollow, by xxx.

\endcomment


As mentioned in [GH, p\,32], there is a weaker equivalence relation than
what we have called isomorphism (positive $A$-module isomorphism with
positive inverse, both of which are isometries). If $M$ and $M'$ are
complete $l^1(\Z)$-modules, we say they are {\it neutrally isomorphic\/}
if there is an $A = l^1(\Z)$-module isometry with isometric inverse from
$M$ to $M'$---in other words, positivity has been left out of the
definition.

Following [GH, p\,32], we say a dimension space (or its sequence of maps),
$M$, is WAT (weakly approximately transitive) if there exists a complete
$l^1(\Z)$-module $N$ given as the completion of the direct limit of the
sequence of maps of the form $\Arrow \times p_i ; l^1 (\Z).l^1(\Z)$ where
$p_i$ are in $l^1(\Z)$ (and can be assumed to be polynomials).

If $\Cal M$ is given as the completion of the (order) direct limit, $\lim \Arrow M_j;
A^{n(j)}.A^{n(j+1)}$ (where $n(j)$ are positive integers, $m_j$ are
$n(j+1) \times n(j)$ matrices with entries in $A^+$ \st after evaluating
all the polynomial entries at $x=1$, the resulting real matrices are
column stochastic), then necessary and sufficient for $\Cal M$ to be WAT is
that there exist a telescoping $t(1) < t(2) < \dots$ \st on defining
$M^{(i)} = M_{t(i+1)-1}\cdot M_{t(i+1)-2}\cdot \dots \cdot M_{t(i)}$,
there exist for each $i$, a row $W_i$ and a column $V_i$ with entries in
$A$, for which
$\sum_i \left\| M^{(i)} - V_i W_i \right \| < \infty$.

In that statement, if we insist that the entries of $V_i$ and $W_i$ belong
to $A^+$, then we have a characterization of AT. There is no requirement
that $W_i V_i$ belong to $A^+$, although this would be desirable.

Now we make a slight weakening of the definition to permit complex
coefficients. We say $M$ is WATC (weakly approximatively transitive over
the complexes), if we permit the entries of $V_i$ and $W_i$ to be in $A
\otimes \C = \C [x^{\pm1}]$ (equivalently, we deal with module
isomorphisms of $M\otimes \C$).

Some invariants for dimension spaces intended for distinguishing
isomorphism classes actually turn out to be neutral isomorphism
invariants---in particular, the mass cancellation-type invariants in [H], under some circumstances, as we will see later.

We will show that if $M$ is the dimension space arising from an ergodic
sequence of hemicirculant matrices, then $M$ is WATC.
If additionally, the matrix size is two, then $M$ is WAT. This
contrasts with the situation for AT, since an example is known (size two)
of an ergodic sequence of circulant matrices which is not AT. (In
particular, WAT does not imply AT.)

\Lem Lemma \onefou. Let $n > 1$ be an integer, and let $0 \leq j_0 < d$ be
integers. Set $T = [j_0, d]\cap \Z$, and let $Z_{kl}$ ($(k,l) \in \brcs{0,1,2,\dots, n-1} \times \brcs{0,1,2,\dots, n-1}$ and $ k\neq l$) be pairwise disjoint subsets of $T$. Then
there exist subsets $U_0, U_1, \dots , U_{n-1}$ of $T$ \st for all $k \neq
l$,
$$
Z_{k,l} \subseteq U_k \cap U_l^c.
$$

\Pf In the unit disk, let $C_0, C_1, \dots, C_{p-1}$ be rays through
the origin with angles $0, \pi/n, 2\pi/n, \dots, (p-1)\pi/n$ respectively.
This yields $n$ open segments partitioning the disk less the union of
the rays. Let $A_{0,+}$ be the open upper half disk, and for each $t$ in $\brcs{0,1,2,\dots, n-1}$,
define $A_{t,+}$ to be $A_{0,+}$ rotated by $t\pi/n$. Then $C_t$ is in
the boundary of $A_{t,+}$. Similarly let $A_{t,-}$ be the complement of
$A_{t,+}$ less $C_t$. Map each $Z_{k,l}$ bijectively into $A_{k,+} \cap
A_{l,-}$ (this is a nonempty sector); we can obviously do this so that the map is bijective on the
union of the $Z_{k,l}$. Call the map $f$, and set $U_t = f^{-1}(A_{t,+})$.
\qed

\Lem Proposition \onefiv. Suppose that $(M_j)$ is an ergodic sequence of hemicirculant matrices \wrt the finite abelian group $H$. Then the corresponding dimension space is WATC. If $H = \Z_2^k$ for some $k$, then the dimension space is WAT.

\Pf Let $n = |H|$ and $Y = \brcs{0,1,\dots, n-1}$. Let $\Arrow F; H. Y$ and $\Arrow E; (Y\times Y) \setminus \Delta. \brcs{1,2,\dots, n^2-n}$ be bijections. Given $\epsilon > 0$, by ergodicity, there exist positive integers $j_0 = d_1 < d_2 < \dots < d_{n^2-n+1} = d$ \st on setting (for unequal $\alpha$ and $ \beta$ in $\hat H$), $Z_{F(\alpha),F(\beta)} = \left[d_{E(F(\alpha),F(\beta))}, d_{E(F(\alpha),F(\beta))+1}\right) \cap \Z$,
$$
\left\|\prod_{j \in Z_{F(\alpha),F(\beta)}} \lambda_{\alpha} (M_j)\lambda_{\beta} (M_j) \right\| < \frac {\epsilon}{n^2}.
$$
By the preceding lemma, there exist subsets $U_{F(\alpha)}$ ($\alpha$ ranging over $\hat H$) of $[j_0,d) \cap \Z$ \st $ Z_{F(\alpha),F(\beta)} \subseteq U_{F(\alpha)} \cap U_{F(\beta)}^c$.

We may write $\prod_{j=j_0}^d M_j = \sum_{\alpha \in \hat H} v_{\alpha} w_{\alpha} \prod_{j=j_0}^d \lambda_{\alpha} (M_j)$. Set
$$
V = \sum_{\hat H} v_{\alpha} \prod_{j \in U_{F(\alpha)}} \lambda_{\alpha} (M_j) \qquad \qquad W = \sum_{\hat H} w_{\alpha} \prod_{j \in U_{F(\alpha)}^c} \lambda_{\alpha} (M_j).
$$
Then
$$\eqalign{
VW & = \sum_{\alpha \in \hat H} v_{\alpha}w_{\alpha} \prod_{j=j_0}^{d} \lambda_{\alpha}(M_j) + \sum_{\alpha \neq \beta} v_{\alpha}w_{\beta} \prod_{j\in U_{F(\alpha)}} \lambda_{\alpha}(M_j) \prod_{j\in U_{F(\beta)}^c} \lambda_{\beta}(M_j)\cr
& = \prod_{j=j_0}^d M_j + \sum_{\alpha \neq \beta} v_{\alpha}w_{\beta} \prod_{j\in U_{F(\alpha)}} \lambda_{\alpha}(M_j) \prod_{j\in U_{F(\beta)}^c} \lambda_{\beta}(M_j),\cr
}$$
Since $\left\| \lambda_{\alpha}(M_j) \right\| \leq 1$ in any event, we have
$$\eqalign{
\left\| \prod_{j\in U_{F(\alpha)}} \lambda_{\alpha}(M_j) \prod_{j\in U_{F(\beta)}^c}\lambda_{\beta}(M_j)\right\| & \leq\left\| \prod_{j \in U_{F(\alpha)}\cap U_{F(\beta)}^c} \lambda_{\alpha} (M_j) \lambda_{\beta} (M_j)\right\| \cr
& \leq \left\| \prod_{j \in Z_{F(\alpha),F(\beta)}} \lambda_{\alpha} (M_j) \lambda_{\beta} (M_j)\right\| < \frac{\epsilon}{n^2}.\cr
}$$
Hence $\| VW - \prod_{j=j_0}^d M_j \| < \epsilon \| \sum v_{\alpha} w_{\alpha} \|/n^2 < \epsilon$.

If $H = \Z_2^k$, then all the entries of each term in the definitions of $V$ and $W$ are real. \qed

As an aside, an easy computation yields $WV = \tr (\prod_{j_0}^M M_j)$,
which shows the former has coefficients in $A^+$. Unfortunately, this
does not (by itself) force $(A_j)$ to be neutrally isomorphic (even over
$\C$) to an AT sequence, because the corresponding sequence is of the form
$W^{(s+1)} V^{(s)}$ (the superscripts indicating the iteration of the
process in the proof of the proposition), and there is no guarantee that
these will have positive (or real) coefficients.

\SecT 2 Hollowness of some circulant sequences

A sequence $(M_j)$, of hemicirculant matrices with entries from $A^+$ and for which the row sums of all the coefficients add to $1$ is {\it hollow\/} if for all $j_0$, and all $\alpha$ in $\hat H \setminus \brcs{\chi_0}$,
$$
\lim_{d \to \infty} \left\|\lambda_{\alpha}\(\prod_{j_0}^{j_0 + d} M_j\)\right\| = 0.
$$
Examples of hollow sequences are easy to construct, as we will show. An obvious example arises when $\prod_{j=j_0}^{j_0 + d} \| \lambda_\alpha (A_j) \| \to 0$ for all $ \alpha \neq \chi_0$, but there are more interesting ones for which this stronger condition fails. For example, if $H = \Z_n$ with $ n \geq 3$, let $P$ denote the standard cyclic permutation matrix of size $n$, and set $M_j = (\slfrac12 (I + x^{2^i}P))^2$; then $\|\lambda_{\alpha} (M_j)\| =1$ for all $\alpha$; however, if we telescope in triples, we find $\|\lambda_{\alpha} (M_{3j} M_{3j+1} M_{3j+2})\| < 1$ when $\alpha \neq \chi_0$, and it easily follows that $(M_j)$ is hollow.

\Lem Lemma \twoone. If $(M_j)$ is a hollow sequence of hemicirculant matrices, then the corresponding dimension space $\lim \Arrow M_j;A^n.A^n$ is AT; there is a telescoping $d(1) < d(2) < \dots$ so that
$$
\sum_l \left\| \prod_{j= d(l)}^{d(l+1)-1} M_j - v_0w_0 \prod_{j= d(l)}^{d(l+1)-1} \lambda_0(M_j )\right\| < \infty.
$$
In particular, the dimension space of $(M_j)$ is isomorphic to that of $(\lambda_0(M_j))$.

\Pf It suffices to show, given $j'$ and $\epsilon>0$, there exists $d$ \st $ \left\| \prod_{j= j'}^{j'+d} M_j - v_0w_0 \prod_{j= j'}^{j'+d} \lambda_0(M_j )\right\| < \epsilon$. For each $\alpha \neq \chi_0$, there exists $d_{\alpha}$ \st $\left\|\prod_{j= j'}^{j'+d} \lambda_{\alpha}(M_j ) \right\| < \epsilon/(n-1)$. If $d= \sup_{\hat H\setminus \brcs{\chi_0}} d_{\alpha}$, then
$$\eqalign{
\left\| \prod_{j= j'}^{j'+d} M_j - v_0w_0 \prod_{j= j'}^{j'+d} \lambda_0(M_j )\right\| & \leq \left\| \sum_{\hat H} \lambda_{\alpha} (\prod_{j= j'}^{j'+d} M_j )v_{\alpha} w_{\alpha}- v_0w_0 \prod_{j= j'}^{j'+d} \lambda_0(M_j )\right\| \cr
& = \left\| \sum_{\hat H\setminus \brcs{\chi_0}} \lambda_{\alpha} \(\prod_{j= j'}^{j'+d} M_j \)v_{\alpha} w_{\alpha}\right\| \cr
& \leq \sum_{\hat H\setminus \brcs{\chi_0}} \frac{ \epsilon}{n-1} \| v_{\alpha} w_{\alpha}\|\cr
& = \epsilon.
}$$
\qed

Now we show that hollowness applies to symmetric and other matrices, with additional constraints. In the symmetric case, the criterion is rather primitive.

The following is reminiscent of Mineka's criterion (for triviality of boundaries), generalizing [GH, Proposition 5.3].

If $p$ and $q$ are Laurent polynomials with real coefficients, define their infimum in the obvious way, $p\wedge q := \sum \min\brcs{(p,x^i),(q,x^i)}x^i$; obviously, if $p$ and $q$ belong to $A^+$, then so does $p\wedge q$, and if $p\wedge q = p$, then $(p,x^i) \leq (q,x^i)$ for all $i$.

\Lem Lemma \twotwo. Suppose that $(M_j = \sum_{g \in H} q_{gj} m_g)$ is an ergodic sequence of hemicirculant matrices \wrt $H$. Sufficient for $(M_j)$ to be hollow is that for all maximal subgroups $K$ of $H$,
$$
\sum_j \sum_{\Set{(g,h)}{gh^{-1} \not \in K}} (q_{gj} \wedge q_{hj})(1) = \infty.
$$

\noindent {\it Remark\/} If $H = \Z_n$ is cyclic (so that the matrices are circulant) of order $n$, it follows that sufficient (but not necessaryunless $n$ is a power of a prime) is divergence of $\sum_j \sum_{i \not \equiv i' \mod n} (q_{ij} \wedge q_{i'j})(1)$.

\Pf Let $\alpha$ be a nontrivial element of $\hat H$, and let $K$ be a maximal subgroup containg the kernel of $\alpha$. Set $a_j = \sum_{\Set{(g,h)}{gh^{-1} \not\in K}} (q_{gj} \wedge q_{hj})(1)$. Since there are fewer than $n^2$ terms in the sum, there exists $(g_j, h_j)$ \st $g_j h_j^{-1} \not \in K$ and $(q_{g_j j} \wedge q_{h_j j})(1) \geq a_j / n^2$. Set $f_j = q_{g_j j} \wedge q_{h_j j}$ and expand
$$\eqalign{
\lambda_{\alpha} (M_j) & = \sum_{g\in H} q_{gj} \alpha (g^{-1}) \cr
& = (\alpha(g_j^{-1}) + \alpha(h_j^{-1}))f_j + (g_{j}- f_j)\alpha(g_j^{-1}) + (h_{j}- f_j)\alpha(h_j^{-1}) + \sum_{g \not \in \brcs{g_j, h_j}} q_{gj} \alpha (g^{-1}); \cr
&\qquad\qquad \qquad\quad \text{as $\sum_H q_{gj} (1) =1$,} \cr
\left\| \lambda_{\alpha}(M_j)\right\| & \leq 1 - (2 -|\alpha(g_j^{-1}) + \alpha(h_j^{-1})|)f_j(1); \quad\text{since $g_jh_j^{-1} \not\in K$ and thus $\alpha(g_jh_j^{-1})\neq 1$,} \cr
& \leq 1 - 2 \sin \frac {\pi}{N} \cdot f_j(1) \leq 1 - a_j \frac{2 \sin \frac {\pi}N}{n^2}\cr
}$$
(Recall that $N$ is the exponent of $H$; we could just as well as replaced it with $n$.) Since $\sum a_j = \infty$, $\prod \lambda_{\alpha} (M_j) \to 0$, verifying hollowness.
\qed

\comment
\Lem Proposition. Suppose that $A_j = \sum_{i=0}^{p-1} q_{ij} B^i$ are circulant matrices of size $p$ with $q_{ij}$ in $S^+$ \st for all $j$, $\sum_i q_{ij}(1) = 1$. Sufficient for the dimension space associated to $(\Arrow A_j; S^p.S^p)$ to be AT is
$$\eqalign{
\sum_j &\sum_{\Cal L} (q_{ij} \wedge q_{p-i,j})(1) = \infty, \quad\text{where } \cr
\Cal L &= \Set{i}{1 \leq i \leq \flo{p/2}; \text{gcd\,}(i,p) = 1;\text{ and (if $p$ is even), $2i$ does not divide a multiple of $p/2$}}.\cr
}$$

\Pf As usual, $\xi = \exp (2\pi \sqrt{-1}/p)$ and $v_k = (1,\xi^k,\xi^{2k}, \dots, \xi^{(p-1)k}\)^T$; the latter is the eigenvector for $B$ with eigenvalue $\xi^k$. If $A = \sum q_i(X)B^i$ (suppressing the $j$, and reminding ourselves that the $q$s are polynomials), then $Av_k = \lambda_k$ where $\lambda_k $ is the complex polynomial $q_0 + \sum \xi^{ki}q_i$. Set $f_i = q_i \wedge q_{p-i}$ so that
$$\eqalign{
\lambda_k &= q_0 + \sum_{i\in \Cal L} (\xi^{ki} + \xi^{-ki})f_i + \sum_{i\in \Cal L} \(\xi^{ki} (q_i - f_i) + \xi^{-ki} (q_{p-i} - f_i)\) + \sum_{i\in \Cal L^c \setminus\brcs{0}} \xi^{ki} q_i \quad \text{so}\cr
\| \lambda_k\| & \leq q_0 (1) + \left\|\sum_{i\in \Cal L} 2 \cos \frac{2k \pi i}{p}f_i \right\|+ \sum_{i\in \Cal L} (q_i - f_i)(1) + \sum_{i\in \Cal L^c \setminus\brcs{0}} q_i(1); \quad\text{ since} \sum q_i (1) = 1, \cr
& \leq 1 - 2\(1-\(1 - \cos\frac{\pi}{p}\)\)\sum_{\Cal L}f_i(1)\cr
& = 1 - 2\cos\frac{\pi}{p}\sum_{\Cal L}f_i(1).
}$$
The inequality from the penultimate line comes from the fact that if $p$ is odd, distance from $|\cos 2j\pi/p|$ to $1$ is minimized at $j = (p-1)/2$, at which the distance is $ 1- \cos \pi/p$ (not $\cos 2\pi/p$ as it would be if $p$
were even).

Restoring the subscript $j$s, we have
$$
\left\| \prod_j \lambda_{kj} \right\| \leq \prod_j \| \lambda_kj\| \to 0
$$
since the hypothesis asserts that $\sum_j f_{ij}(1) = \infty$.

Thus for all $\epsilon > 0$ and any $j_0$, there exists $M_0 \geq j_0$ \st for all $k \neq 0$, $\left\|\(\prod_{j_0}^{M_0}A_j\) v_k\right\| < \epsilon $. We may write the standard basis $p e_l = \sum_{k \in \Z_p} \xi^{-lk}v_k = v_0 + \sum_{k=1}^{p-1} \xi^{-lk}v_k$, so that for each $l$,
$$
\left\|\(\prod_{j_0}^{M_0}A_j\) e_l - \frac1p\(\prod_{j_0}^{M_0}A_j\) v_0\right\| \leq \frac{p-1}p
\epsilon.$$
In particular, any two columns of $\prod_{j_0}^{M_0}A_j$ differ from each other by a vector of norm less than $2\epsilon$. This yields an approximate factorization of $\prod_{j_0}^{M_0}A_j$ as $\prod Q_j \cdot v_0 v_0^T/p$, where $Q_j = \sum_i q_{ij}$ (and $v_0 v_0^T$ is the matrix all of whose entries are $1$) with error at most $2p\epsilon$. Since we can adjust $\epsilon $ as we proceed, we can ensure that the sum of the errors over all subsequent telescopings converge. \qed
\endcomment

The matrix $M = \sum_{g \in H} q_g m_g$ is symmetric if and only if for all $g$, $q_g = q_{g^{-1}}$.

\Lem Proposition \twothr. Let $(M_j = \sum_{g\in H} q_{gj} m_g)$ be an ergodic sequence of hemicirculant matrices \wrt to $H$. Suppose in addition that each $M_j$ is symmetric. Then $(M_j)$ is hollow if for all maximal proper subgroups $K$ of $H$,
$$
\sum_j \sum_{\Set{g \in H}{g^2 \not \in K}} q_{gj} (1) = \infty. \tag \dag
$$
In particular, if $|H|$ is odd and the $M_j$ are symmetric, then $(M_j)$ is hollow.

\Pf By symmetry, $q_{gj} \wedge q_{g^{-1}j} = q_{qj}$ and $g(g^{-1})^{-1} = g^2$. Thus
$$
\sum_j \sum_{\Set{g,h}{gh^{-1} \not \in K}} (q_{gj} \wedge g_{hj})(1) \geq \sum_j \sum_{\Set{g}{g^2 \not \in K}} q_{gj}(1),
$$
and the previous lemma applies.

Let $f_j$ in $H$ be such that $q_{f_j j}(1) $ is maximal among $\brcs{q_{gj}(1)}$. Then $f_j M_{j} = \sum_{g \in H} q_{f^{-1}g} m_g$, and if $f_j^2 = 1$, then $f_j M_j$ is still symmetric but now the maximal coefficient in $f_j M_j (1) = \sum q_{f_jg}(1) m_g$ occurs at the identity coefficient. If $f_j^2 \neq 1$ (as will be the case when $|H| $ is odd and $f_j$ is not the identity), then the coefficients of $M_j(1)$ have maxima at two distinct points, $q_{f_j j} (1)$ and $q_{f_j j}(1)$. This together with ergodicity is more than enough to guarantee that the condition in Lemma \oneone(iii),
$$
\sum_{j} \sum_{g\not\in K} q_{gj}(1) = \infty \tag *
$$
for all maximal subgroups $K$ (in place of the maximum coefficient occuring at the identity, which is assumed there).

Now suppose $|H|$ is odd. Then $H/K$ has odd order, so $g^2\in K$ entails $g \in K$. Thus ($\dag$) is equivalent to ($*$) (for all maximal subgroups $K$), hence $(M_j)$ is hollow.
\qed

The condition, $g^2 \not \in K$ (as opposed to $g \not \in K$), really is significant. The standard example of a non-AT but AT(2) ergodic action [GH] is given by size two circulant matrices ($H = \Z_2$)---necessarily symmetric---$M_j =\frac 12 \(\smallmatrix 1 & x^{5^{j}}\\ x^{5^{j}} & 1 \\ \endsmallmatrix \)$. The sequence $(M_j)$ is not hollow, as it is not even AT.

In many cases, $(M_j^2)$ is hollow; this occurs if (for example) $H = \Z_n$ (cyclic)
and $M_j = \frac 12 (1 + x^{2^j}m_{[1]})$ where $[1]$ is the generating element of $
Z_n$ (so $m_{[1]}$ is just the usual cyclic permutation matrix of size $n$). Here $(M_j)$ is not hollow, as is easy to check. Less easy to check is that
$\| \prod_{j=0}^d (1 + \xi^l x^{2^j})^2\| = \oh{4^d}$ if $l $ is not divisible by $n$ (in fact, convergence is faster than $\oh{s^d}$ for some number $s < 4$, possibly as low as $s=2$), so $(M_j^2)$ is hollow. However, there exist sequences $(M_j) $ for which hollowness of powers does not occur.

For example, suppose $H = \Z_n$ again, let $k$ be a positive even integer exceeding two. Set $M_j = \frac 12 (1 + x^{k^j} m_{[1]})$. Then $(M^{s}_j)$ is not hollow for $s < k$, but is hollow for $s \geq k$. More drastically, if $M_j = \frac 12 (1 + x^{j!} m_{[1]})$ (this example has extravagantly large gaps in the exponents), then $(M_j^s)$ is not hollow for any $s$. More drastically still, if $M = \frac12 \(\smallmatrix 1 + x^2 & x\\ x & 1+x^2 \\ \endsmallmatrix\)$ and $M_j = M(x^{3^j})$, and $\Arrow f; \N.\N$ is any function whatsoever, then $\(M_j^{f(j)}\)$ is not hollow. (This is easily deduced from the signs of the coefficients in the products $\lambda_1\(\prod M_j^{f(j)}\) = \prod\frac12 (1-x^{3^j})^{f(j)}$: the coefficient of $x^k$ is $(-1)^k$ times the coefficient of $x^k$ in $\lambda_0\(\prod M_j^{f(j)}\) = \prod\frac12 (1+x^{3^j})^{f(j)}$.)

The following is practically tautological.

\Lem Proposition \twofou. Let $(M_j)$ an ergodic $H$-hemicirculant sequence. The following are equivalent.
{\par}
\item{(i)} $(M_j)$ is hollow;
\item{(ii)} the kernel of the natural map from the dimension space of $(M_j)$ to that of $(\lambda_0 (M_j))$ is zero;
\item{(iii)} for all $g $ in $H$, the automorphism of the dimension space $(M_j)$ given by $[f,k] \mapsto [m_g f,k]$, is the identity.

Suppose $\Arrow g; \N.\N$ is a function, and set $M_j = \frac12 \(\smallmatrix 1 & x^{g(j)}\\ x^{g(j)} & 1 \\ \endsmallmatrix\)$. For example, $g(j)$ could be the $j$th Fibonacci number. Under some conditions (including this example), the sequence $(M_j)$ is hollow; under somewhat weaker conditions, $(M_j)$ is AT but not necessarily hollow. This requires several lemmas. First, we have a special case in which the method used to show $(M_j^2)$ is AT sometimes works to show $(M_j)$ is AT.

\Lem Lemma \twofiv. Let $(M_j)$ be an ergodic sequence of $H$-hemicirculant matrices. Suppose that for all $\epsilon > 0$, and all sufficiently large integers $k$, there exists a positive integer $ d\equiv d(\epsilon, k)$ together with a subset $R \subseteq T:= \brcs{k,k+1, \dots, k+d}$ \st for all $\alpha \neq \beta$ in $\hat H$,
$$
\left\| \prod_{j\in R} \lambda_{\alpha} (M_j) \prod_{j\in T \setminus R} \lambda_{\beta} (M_j) \right\| < \epsilon.
$$
Then $(M_j)$ is AT.

\Pf Set $M = \prod_T M_j$. Set $V = \sum_{\alpha\in \hat H} v_{\alpha}\prod_{j\in R} \lambda_{\alpha} (M_{j})$ and $W = \sum_{\alpha\in \hat H} w_{\alpha}\prod_{j \in T\setminus R} \lambda_{\alpha} (M_{j})$. As in the argument of Theorem \onethr, the entries of both $V$ and $W$ are in $A^+$ (explicitly, $V = \sum q_g e_g$ where $\prod_R M_j = \sum q_g m_g$, and $W = \sum q_g' e_g$ where $\prod_{T\setminus R} M_j = \sum q_g' m_g$). Also as in that argument,
$$\eqalign{
VW & = \sum_{\alpha \in \hat H} v_{\alpha} w_{\alpha} \lambda_{\alpha} (M) + \sum_{\alpha \neq \beta} v_{\alpha}w_{\beta}\prod_{j \in R} \lambda_{\alpha}(M_j)\prod_{j\in T \setminus R} \lambda_{\beta} (M_j); \quad \text{thus,} \cr
\left\| VW - M \right\| & \leq |H| (|H|-1) \epsilon.\cr
}$$
This yields suitable approximate factorizations of the products of the $M_j$, so $(M_j)$ is AT. \qed

\Lem Lemma \twosix. Let $\Arrow g; \N.\N$ be a function satisfying the following conditions:
{\par}
There exists $S \subseteq \N$ together with $\epsilon_i (s) \in \brcs{0, \pm 1}$ (for $s\in S$), \st {\par}
\item{(i)} for all $s \in S$, we may write $g(s) = \sum_{i=1}^{N-1} \epsilon_i (s) g(s-i)$ where $\sum_i \epsilon_i(s) $ is even;
\item{(ii)} on setting $b(s) = 1 + \sum |\epsilon_i (s)|$, we have $\sum 2^{-b(s)} = \infty$;
\item{(iii)} on defining $\supp g(s) = \brcs{s} \cup \Set{s-i}{\epsilon_i (s) \neq 0}$, we have $ s\neq s' \in S$ entails $\supp g(s) \cap \supp g(s') = \emptyset$.
{\par}
\noindent Then for all integers $k$,
$$
\lim_{N\to \infty}\left\| \prod_{j=k+1}^{N+k} \frac{1-x^{g(j)}}2\right\| = 0.
$$

\Pf Define $S_k = \Set{s}{\supp g(s) \cap \brcs{1,2,3,\dots, k} = \emptyset}$. Obviously $\sum_{s \in S_k} 2^{-b(s)} = \infty$. For $\epsilon > 0$, there exists $N\equiv N(\epsilon)$ sufficiently large that if $A_N = \Set{s\in S_k}{\supp g(s) \subseteq \brcs{k+1,\dots, N+k}}$, then
$$
\prod_{s \in A_N} \(1 - \frac 1{2^{b(s)}} \) < \epsilon.
$$
Define $V_{N} = \prod_{k+1}^{N+k} \brcs{0,1}$. Say that $v$ in $V_N$ {\it runs through $\supp g(s)$} if
$$
v(i) = \cases 1& \text{$\epsilon_i (s) = 1$} \\
0& \text{$\epsilon_i (s) = -1$}. \\
\endcases
$$
Let
$${\Cal S} = \Set {v\in V_N}{\text{$v$ runs through at least one $ g(s)$ with $s\in A_N$.} }
$$
Then $|V_N \setminus {\Cal S}| < \epsilon |V_N| = 2^N \epsilon$, so $|\Cal S| \geq 2^N (1-\epsilon)$.

Define $\Arrow \phi;{\Cal S}. V_N$ as follows. For $v \in {\Cal S}$, find the smallest $s$, denoted $s_0$, \st $s \in A_N$ and $v$ runs through $\supp g(s)$; define
$$
\phi(v)(i) = \cases
1-v(i) & \text{if $i \in \supp g(s_0)$}\\
v(i) & \text{else.}\\
\endcases
$$
Then $\phi$ is one to one, and
$$\eqalign{
\sum_{k+1}^{N+k} \phi(v)(i) & \equiv 1 + \sum v(i) \mod 2\cr
\sum_{k+1}^{N+k} \phi(v)(i) g(i)& = \sum v(i)g(i).\cr
}$$
Now
$$
\left\|\prod (1- x^{g(j)}) \right\| = \sum_{j=k+1}^{k+N} \left| \sum_{ \sum v(i)g(i) = j \& v \in A_N} (-1)^{\sum v(i)}\right|.
$$

If $v \in \Cal S$, then $\phi(v)$ and $v$ have opposite parities (that is, $(-1)^{\sum v(i)} + (-1)^{\sum \phi(v)(i)} = 0$). This means that the contribution of such a $v$ to the total mass is zero. Hence
$$
\left\| \prod \( 1 - x^{g(j)} \)\right\| \leq |V_N \setminus {\Cal S}| < 2^N \epsilon.
$$
\qed

Part of the hypotheses above practically require $g(n) = \oh{2^n}$ (since we only use the hypothesis $g(s) = \sum_{i < s} \epsilon_s(i) g(i)$ for $s \in S$, there is no general requirement about the growth of $g$ on all of $\N$; however, in examples, it is easier to assume the recurrence relations hold for all $n$). This aspect cannot be much improved, since (for example), if $g(n) > \sum_{i < n} g(i)$ for all $n$ (e.g., $g$ increasing and $g(n) \geq 2^n$), then no cancellation occurs, and thus $\| \prod (1-x^{g(n)})/2 \| = 1$ (put another way, if $m =
\sum_{n} \eta(n) g(n)$ where $\eta (n) \in \brcs{0,1}$, then $\eta$ is uniquely determined by $m$).

The parity hypothesis (condition (i)) is crucial; for example, if each $g(n)= g(n-1) + g(n-2) + g(n-3)$ for $n > 3$ and each of $g(1), g(2), g(3) $ is odd, then all $g(n) $ are odd, hence for any product of distinct terms of the form $(1 - x^{g(j)})/2$, the norm is $1$, since no cancellation occurs; in particular, $(M_j)$ is not hollow in this case.

On the other hand, if $g$ satisfies $g(j) = g(j-1) + g(j-2)$ (or anything else with an even number of terms), then we can take $S$ to be $3\N$ (to avoid overlaps in the supports), and the hypotheses are all satisfied.

This yields a proof that the $\Z$-action obtained from the sequence $\brcs{\frac 12 \(\smallmatrix 1 & x^{g(i)} \\ x^{g(i)} & 1 \\ \endsmallmatrix\)}$ is hollow. In the case that $g(j)$ is the $j$th Fibonacci number, $(M_j)$ is isomorphic to the AT system $((1+x^{g(j)})/2)$. We now drop the parity condition on $\sum_i \epsilon_i(s)$.

\Lem Proposition \twoeig. Let $\Arrow g; \N . \N$ be a function satisfying conditions (i)--(iii): {\par}
There exists $S \subseteq \N$ together with $\epsilon_i (s) \in \brcs{0, \pm 1}$ (for $s\in S$), \st {\par}
\item{(i)} for all $s \in S$, we may write $g(s) = \sum_{i=1}^{N-1} \epsilon_i (s) g(s-i)$;
\item{(ii)} on setting $b(s) = 1 + \sum |\epsilon_i (s)|$, we have $\sum 2^{-b(s)} = \infty$;
\item{(iii)} on defining $\supp g(s) = \brcs{s} \cup \Set{s-i}{\epsilon_i (s) \neq 0}$, we have $ s\neq s' \in S$ entails $\supp g(s) \cap \supp g(s') = \emptyset$.
{\par}
\noindent Then {\par}
\noindent(a) $\(M_j := \frac 12 \(\smallmatrix 1 & x^{g(i)} \\ x^{g(i)} & 1 \\ \endsmallmatrix\)\)$ is AT.
{\par}
\noindent (b) $\(M_j\)$ is hollow if additionally, $\sum_i \epsilon_i (s)$ is even for all $s \in S$. In this case, the dimension spaces corresponding to $(M_j)$ and to
$((1 + x^{g(j)})/2)$ are isomorphic.

\Pf (a) Suppose $s$ belongs to $S$; we have a relation of the form $g(s) + \sum_{i \in E_- (s)} g(s-i) = \sum_{i \in E_+ (s)} g(s-i)$, where $E_+(s) = \Set{i \in \N}{\epsilon_i (s) = 1}$ and $E_-(s) = \Set{i \in \N}{\epsilon_i (s) = -1}$. Since $\epsilon_i (s) \neq 0$ entails $i < s$, each set is finite. Since $g$ is positive-valued, $E_+ (s)$ is nonempty.

This equality yields that at least two terms cancel in each of $(1+x^{g(s)})\prod_{E_-(s)} (1+x^{g(s-i)}) \prod_{E_+(s)} (1-x^{g(s-i)})$ and $(1-x^{g(s)})\prod_{E_-(s)} (1-x^{g(s-i)}) \prod_{E_+(s)} (1+x^{g(s-i)})$, hence
$$
\left\| \frac{(1\pm x^{g(s)})}2 \prod_{i \in E_-(s)} \frac{(1\pm x^{g(s-i)})}{2}\prod_{i \in E_+(s)} \frac{(1\mp x^{g(s-i)})}{2}\right\| \leq 1 - \frac 2{2^{b(s)}}.
$$
Since $\lambda_{0} (M_j) = (1+ x^{g(j)})/2$ and $\lambda_1 (M_j) = (1-x^{g(j)})/2$, we will verify the conditions of Lemma \twofiv\ if we select $d$ large enough that we can choose $R$ inside $T$ so that both $R$ and $T\setminus R$ contain enough terms of $S$ to guarantee that the two products are arbitrarily small; this is possible since $\prod (1- 2^{-b(s)}) = 0$, as a consequence of $\sum 2^{-b(s)} = \infty$.

\noindent (b) Since the matrices are two by two and circulant, the only relevant eigenvalue of $M_j$ is $(1- x^{g(j)})/2$; so hollowness follows from Lemma \twosix{} (and does not depend on the argument in (a)).
\qed

For example, if $g(j)$ is the $j$th Fibonacci or Lucas number, then $(M_j)$ is isomorphic to the AT system $((1+x^{g(j)})/2)$, since in either case, $g(j) = g(j-1) + g(j-2)$, and we can set $S = 3\N$. On the other hand, if $g(j) = g(j-1) + g(j-2) + g(j-3)$ for all sufficiently large $j$, set $ S = 4\N$; (a) yields that the resulting sequence $(M_j)$ is at least AT, but need not be hollow, as we observed above.

Roughly speaking, if $g(j)$ is $\oh{2^j}$ (and satisfies a recurrence relation as above), then the resulting $\(M_j = \frac 12 \(\smallmatrix 1 & x^{g(j)}\\x^{g(j)} & 1 \\ \endsmallmatrix\)\)$ tends to be at least AT, and could be hollow. A limiting case (and not covered by the hypotheses) occurs with $g(j) = 2^j$, for which ATness of the resulting sequence is due to [DQ]. It is known that if $g(j ) = 5^j$, then the sequence is not AT [GH].
\comment
WRONG	!
\Lem Proposition. Suppose that $(M_j)$ is an ergodic sequence of $H$-hemicirculant matrices which is isomorphic to an AT sequence with no outer automorphisms of finite order. Then $(M_j)$ is hollow.

\Pf For each $g$ in $H$, $m_g$ induces an automorphism of the dimension space,
via $e_g[f, k] = [e_g \cdot f,k]$ (this commutes with the corresponding $M_k$). By hypothesis, the automorphism must be trivial, that is, in particular, $[e_g - 1,k] = 0$, and hollowness easily follows.
\qed

Examples include odometers (corresponding to supernatural numbers, probably all odometers).
\endcomment

\def\Tr{\text{Tr\,}}
\let\Tr=\tr

\Lem Proposition \twonin. If $(M_j) $ is an ergodic sequence of hemicirculant matrices, then $(\lambda_0 (M_j) M_j ) = (\lambda_0(M_j)) \otimes (M_j)$ is hollow, and thus
$(\lambda_0 (M_j) \cdot M_j ) \iso (\lambda_0(M_j)^2)$.

\Pf We note for $\alpha$ in $\hat H$, that $\lambda_{\alpha}(\lambda_0 (M_j)M_j) = \lambda_{\alpha}(M_j)\lambda_0 (M_j)$; thus with $\beta =
\chi_0$ and $\alpha \neq \chi_0$, setting $N_j = \lambda_0 (M_j) M_j $, we have
$$
\lambda_{\alpha} \(\prod_{j_0}^{j_0 + d} N_j\) = \prod_{j_0}^{j_0 + d} \lambda_{\alpha}(M_j)\lambda_{\beta} (M_j),
$$
which goes to zero as $d\to \infty$, by Lemma \onetwo.
\qed

If instead, we ask about multiplication by suitable values of the trace, an extra condition is required. For size two matrices, there is an analogous result even in the noncommutative case, arising from $N^2 = N\tr N - \I \det N$; by suitably telescoping, we can arrange that the determinant terms are as small as we like. As a consequence, if $(N_j)$ is an ergodic sequence of $2\times 2$ matrices, then there exists a telescoping so that $((N^{(i)})^2) \iso (\tr(N^{(i)})\cdot N^{(i)})$. This is not as useful as it might appear---there is no assumption that the $N_j$ commute, so the isomorphism class of $((N^{(i)})^2)$ might depend on the choice of telescoping.

\Lem Proposition \twoten. Suppose that $(M_j)$ is an ergodic sequence of hemicirculant matrices. Then $(M_j^2)$ is hollow if and only if there exists a telescoping $0 = n(1) < n(2) < \dots$ with $M^{(i)} = \prod_{j=n(i)}^{n(i+1)-1} M_j$ \st the sequence $(\Tr(M^{(i)})\cdot M^{(i)})$ is hollow.

\Pf 
For $\alpha $ in $\hat H \setminus{\chi_0}$ and $J$ a finite subset of $\N$, we observe that $\Tr(\prod_{j \in J} M_j) = \sum_{\gamma \in \hat H} \prod_J \lambda_{\gamma} (M_j)$. Thus
$$\eqalign{
\lambda_{\alpha} \(\Tr \(\prod_{j \in J} M_j \)\cdot\prod_{j\in J} M_j \) & = \sum_{\gamma \in \hat H} \prod_J \lambda_{\gamma} (M_j)\lambda_{\alpha} (M_j)\cr
& = \prod_J \lambda_{\alpha} (M_j^2) + \sum_{\alpha \neq \gamma} \prod_{j\in J} \lambda_{\gamma} (M_j)\lambda_{\alpha} (M_j). \cr
}$$

Assume $(M_j^2)$ is hollow. Given any $j = j_0$ we can make the first term on the right side as small as we like by increasing the $d$ in $J = \brcs{j_0,j_1, \dots, j_0 + d}$; by \onetwo, and similarly increasing the $d$ (further if necessary), the second summand can be made as small as we like, hence the desired telescoping exists.

Conversely, assume $(\Tr(M^{(i)})\cdot M^{(i)})$ is hollow. By a further telescoping (which amounts to deleting some of the $n(i)$s and re-indexing), we may assume that the second summand and the entire left side are arbitrarily small for $J = \brcs{n(i),n(i)+1, \dots, n(i+1)-1}$ (and sufficiently large $i$). Thus the first summand must be small, which yields hollowness of $(M_j^2)$.
\qed

\comment
for tensor product

$e_g \otimes e_h \mapsto e_{gh}$

corresponding map on eigenvectors given by
$$
v_{\alpha} \otimes v_{\beta} \mapsto \cases 0 &\text{if $\alpha \neq \beta$}\\
v_{\alpha} & \text{if $\alpha = \beta$}}
\endcases
$$

kernel of downmap after extension\\ by complex scalars is span v \otimes v differing alphas and betas, makes it easier to show map vertical map is an iso (inverse map etc)

discussion of connection to finite abelian group actions, topological then measure-theoretic setting, actions of H on dimension space (nontrivial) corresponding to dual action on vN algebras, hollow iff induced action is trivial, , example with Cosh,
nonisomorphism via invariant, and with powers not being hollow,
\endcomment

\SecT 3 Tensor products

\noindent Let $M$ and $N$ be hemicirculant matrices \wrt the same group $H$. We may form $M \otimes N$, which is now hemicirculant \wrt $H \times H$. If $(M_j)$ and $(N_j)$ are ergodic sequences of hemicirculant matrices, we may form $(M_j \otimes N_j)$, which corresponds to the tensor product of the corresponding dimension spaces.

There are natural $A$-module maps in this situation. Define $\Arrow \theta ;A H \otimes AH . A H$ via $\theta (e_{g} \otimes e_h) = e_{gh}$; this is a positive, onto $A$-module map \st $\theta \circ (M \otimes N) = MN \circ \theta $ (here $M$ and $N$ are $H$-hemicirculant matrices); this requires a routine verification. We may complexify $A_{\C} = A \otimes \C = \C[x,x^{-1}]$ as usual, and extend $\Arrow \theta_{\C} ;A_{\C} H \otimes A_{\C}H . A_{\C} H$. It is easy to verify from the definitions that for $\alpha$ and $\beta$ in $\hat H$,
$$\theta_{\C} (v_{\alpha} \otimes v_{\beta}) = \cases
v_{\alpha} & \text{if $\alpha = \beta$}\\
0 & \text{if $\alpha \neq \beta$.}\\
\endcases
$$
Since $\theta_{\C}$ is actually defined over the complexes, it easily follows (since $\brcs{v_{\alpha}}_{\alpha \in \hat H} $ is a basis for $\C H$), that $\ker \theta_{\C}$ is the free $A_{\C}$-module on $\brcs{v_{\alpha} \otimes v_{\beta}}_{\alpha \neq \beta}$.

The following is a fragment of the commuting diagram obtained from the map $\theta$; the splitting map $\psi$ is given by $\psi (e_g) = \frac1{|H|}\sum_{h \in H} e_h \otimes e_{gh^{-1}}$ (the corresponding complexification of $\psi$ sends $v_{\alpha}$ to $v_{\alpha}\otimes v_{\alpha}$).

$$
\diagram
A^{n^2}&\rTo^{M_{j-1}\otimes N_{j-1}} & A^{n^2} & \rTo^{M_{j}\otimes N_{j}} &
A^{n^2} & \rTo^{M_{j+1}\otimes N_{j+1}} \\
\dTo^{\theta}\uTo_{\psi} &&\dTo^{\theta}\uTo_{\psi} & & \dTo^{\theta}\uTo_{\psi} & \\
A^n&\rTo^{M_{j-1}N_{j-1}} & A^{n}& \rTo^{M_{j}N_{j}} & A^{n}&
\rTo^{M_{j+1}N_{j+1}} \\
\enddiagram
$${\vskip 3pt}
\noindent Now $\theta$ induces an $A$-module map from the dimension space of $(M_j \otimes N_j)$ to that of $(M_j N_j)$, obviously of the form $\Theta\: [z, k] \mapsto [\theta (z),k]$; this defines $\Theta$ on the algebraic direct limits; it clearly extends to the $l^1$ completions, that is, the dimension spaces themselves. It is clearly of norm $1$, positive, onto on positive cones, and sends the unique invariant measure on the domain to that on the image; it also has a splitting, which we will use shortly.

Under some circumstances, $\Theta$ will be an isomorphism. Suppose $z$ in $AH$ is sent to zero by $\theta$; we wish to determine whether for all $k$, and for all $\epsilon > 0$, there exists $d $ \st $\left\| \prod_{j= k}^{k+d} (M_j \otimes N_j)z\right\| < \epsilon$. We can write $z $ as an ($A_{\C}$) linear combination of terms of the form $v_{\alpha} \otimes v_{\beta}$ with $\alpha \neq \beta$. Next, we note that $\(\prod_{j= k}^{k+d} (M_j \otimes N_j)\)v_{\alpha} \otimes v_{\beta} = \(\prod_{j= k}^{k+d} \lambda_{\alpha}(M_j) \lambda_{\beta}(N_j)\)v_{\alpha} \otimes v_{\beta}$. Define the following property,
$$
\text{For all $\alpha \neq \beta \in \hat H$, for all $k$, }\quad\lim_{d \to \infty} \prod_{j= k}^{k+d}\lambda_{\alpha}(M_j) \lambda_{\beta}(N_j) = 0. \tag 1
$$
When (1) holds, it follows that the elements $[z,k]$ in the dimension space of $(M_j \otimes N_j)$ are all zero.

Now assume (1) and suppose $[y,k]$ is such that $\|\theta y\| < \epsilon$. In $A_{\C} H$, we can write $y = z + x$ where $y \mapsto z$ is the projection of $A_{\C} H$ onto $\ker \theta$ (this has norm $1$, but all we need is that it has finite norm, which is trivial here), and there exists a constant $K$
(independent of $y$) \st $\| x \| \leq K\epsilon$. Then it easily follows that the image of $[y,k]$ in the dimension space of $(M_j \otimes N_j)$ has norm at most $K \epsilon$. We can avoid this last argument if we employ the splitting map, $\Arrow \psi; AH . AH \otimes AH$. Then $\psi$ induces a map $\Psi$ from the dimension space of $(M_jN_j)$ to that of $(M_j \otimes N_j)$ and $\Theta \circ \Psi$
is the identity on the dimension space of $(M_jN_j)$.
Hence we obtain, using Lemma \onetwo:

\Lem Lemma \throne. Let $(M_j)$ and $(N_j)$ be ergodic sequences of $H$-hemicirculant matrices. The map $\Theta$ from the dimension space of $(M_j \otimes N_j)$ to that of $(M_jN_j)$ is an isomorphism of dimension spaces if (1) holds. In particular, this occurs if $M_j = N_j$---that is, $(M_j^2)$ and $(M_j \otimes M_j)$ yield isomorphic dimension spaces, necessarily AT.

Another consequence of the arguments used before: if (1) holds, then $(M_j N_j)$ is AT (simply use the argument of \twosix, constructing $V_j$ from $M_j$ and $W$ from $N_j$). However, when $M_j$ and $N_j$ are quite different from each other, it may be difficult to decide whether (1) holds. For example, in general (1) will fail if $N_j = M_j^T$. This does not imply $(M_j \otimes M_j^T) \not\iso (M_jM_j^T)$, simply that the method of proof fails. However, after we remind the reader of an invariant developed in [H, Section 6], we will give a large class of examples for which $(M_j \otimes M_j^T) \not\iso (M_jM_j^T)$.

\SecT 4 Isomorphism and nonisomorphism

Let $(M_j = \sum_{g\in H} q_{gj}m_g) $ be an (ergodic) sequence of H-hemicirculant matrices, with $v_0 = (|H|)^{-1}(1 \ 1 \ \dots \ 1)^T$ and $w_0 = |H| v_0^T$ as usual. The row $w_0$ induces (after a relabelling) the map $\Arrow \rho; AH. A$ sending $\sum_{g \in H} q_g m_g $ to $\sum q_g$. Then we obtain a commuting diagram (not just approximately commuting)
$$\diagram
A^{n}&\rTo^{M_{j-1}} & A^{n} & \rTo^{M_{j}} &
A^{n} & \rTo^{M_{j+1}} \\
\dTo^{\rho}\uTo_{\sigma} &&\dTo^{\rho}\uTo_{\sigma} & & \dTo^{\rho}\uTo_{\sigma} & \\
A&\rTo^{\times \lambda_0(M_{j-1})} & A& \rTo^{\times \lambda_0(M_{j})} & A&
\rTo^{\times \lambda_0(M_{j+1})} \\
\enddiagram \tag \dag
$$
where the horizontal maps on the second row are multiplications by the indicated elements of $A^+$, in this case, the large eigenvalue of the $M_j$, $\lambda_0 (M_j) = \sum q_{gj}$. Then $\rho$ extends to $P$ (capital $\rho$), a positive $A$-module map from the dimension space of $(M_j)$ to the AT dimension space of $(\lambda_0 (M_j))$. With the same elementary techniques as those of the earlier tensor product constructions, it follows that $P$ is an isomorphism if and only if $(M_j)$ is hollow. [This is made easier by the presence of the splitting map for $P$, induced by $\Arrow \sigma; A. AH$ given by $\sigma (c) = \frac 1{|H|}c\sum_{g\in H} e_g$.]

When $(M_j)$ is hollow, it also is immediate that for all $j_0 \geq 0$,
$$
\lim_{d \to \infty} \left\| \tr \(\prod_{j=j_0}^{j_0 + d} M_j\) - \lambda_0 \(\prod_{j=j_0}^{j_0 + d} M_j\) \right\| = 0.
$$
(The norm is the usual $l^1$ norm---aka the total variation norm---on the coefficients.) Since $\lambda_0$ is multiplicative on $H$-hemicirculant matrices, this yields that the dimension space associated to $(M_j)$ is isomorphic (in all senses) to that of a telescoped sequence $\( \tr \(\prod_{j=n(i)}^{n(i+1)-1} M_j\) \)_i$, where $n(i)$ is a suitable strictly increasing sequence of positive integers.

We can ask, if $(M_j)$ is not hollow, then can it be isomorphic to the AT system $(\lambda_o (M_j))$? In many cases (and possibly all), it is relatively easy to give a negative answer. This uses an invariant (actually, a class of invariants) introduced in [H], related to mass cancellation.

In [H], it was intended to be used to distinguish pairs of AT systems, but in fact, can be extended to any dimension space.
Let $\Arrow N_j; A^{n(j)}. A^{n(j+1)}$ be a sequence of matrices with entries from $A^+$ \st all column sums of $ N_j(1)$ are $1$. Let $\brcs{p_i}_{i\in I}$ (in examples, $ I = \N$, but could be any infinite set) be a collection of elements of $A_{\C} = \C[x,x^{-1}]$ (not $A^+$) each with norm $1$. Associate to the sequences $\brcs{p_i}, (N_j)$ a real number $s \equiv s\(\brcs{p_i}, (N_j)\)$ in $[0,1]$ as follows.

Fix a positive integers $l$ and $d$, and $p \in \brcs{p_i}$, and define the operator $l^1$ norm
on the products $p N_{l+d} \cdot N_{l+d-1} \cdot \dots \cdot N_l$
(using the norm inherited from the limiting dimension space). Define
$$
s_{l,p} = \lim_{d \to \infty} \tripnorm pN_{l+d} \cdot N_{l+d-1} \cdot \dots \cdot N_l xxx.
$$
Since each $N_j$ has norm $1$, it easily follows that the limit exists and is a number in the unit interval. Now set $s_l = \inf_{p \in \brcs{p_i}} s_{l,p}$. As $l \mapsto s_l$ is increasing, $s\(\brcs{p_i}, (N_j)\)):= \lim_{l \to \infty} s_l$ is well-defined and in the unit interval. The invariant is $s\(\brcs{p_i},\cdot\)$. We have to check that it is an invariant; in fact, it is even an invariant for neutral isomorphism as well.

First, we note that it is invariant under telescoping and deletion of a finite set of $N_j$. Now suppose we have an approximately commuting diagram as below (here the $M_j$ need not be hemicirculant, or even square), meaning, \wrt the $l^1$ norms induced the invariant linear functionals, we have $\sum \tripnorm S_j R_j - M_j xxx < \infty$ and $\sum \tripnorm R_{j +1} S_j - N_j xxx < \infty$. (By repeatedly telescoping, we may arrange this out of any approximately commuting isomorphism.)

Given $\epsilon > 0$, there exists $l$ \st both sums $\sum_{j \geq l} \tripnorm S_j R_j - M_j xxx, \sum_{j\geq l} \tripnorm R_{j +1} S_j - N_j xxx < \epsilon$ (the norms are on different spaces, but the notation would become cumbersome if we distinguished them). Then to within $\epsilon$, for any polynomial $p$, $p M_{l+d}M_{l+d-1}\dots M_d \sim p S_{l+d} N_{l+d-1}\dots N_l R_l$. Since (we can assume) $R_j$ and $S_j$ have respective norms at most $1$, it follows $s_{l,p} (\brcs{M_j}) \leq \epsilon + s_{l,p} (\brcs{N_j})$. Hence the corresponding result holds dropping the $l$ and taking infima over all the elements $p $ of $\brcs{p_i}$, and of course, then we interchange the roles of $M_j$ and $N_j$.

As an aside, we almost always need an infinite choice of $\brcs{p_i}$. For example, the AT system given by $\(\frac12 (1+ x^{2^j})\)$ represents the dyadic odometer. If $\brcs{p_i} $ is just the singleton $\brcs{(1-x)/2}$, then $s_0 = 0$ (almost complete mass cancellation) but $s_l = 1$ (no mass cancellation) for $l \geq 1$, hence $s=1$, which in this case is uninformative. On the other hand, if
$p_i = (1-x^{2^{i-1}})/2$, then $s_l = 0$ for all $l$ and thus $s = 0$, which is frequently useful.

For example, it allows us to distinguish the sequence $(M_j = \frac12 (I+ x^{2^j}P))$ (where $P$ is the cyclic permutation matrix of size $n$, corresponding to $ H = \Z_n$) from the dyadic odometer. We note that $\lambda_0 (M_j) = \frac12(1+ x^{2^j})$, so the sequence in the bottom line of ($\dag$) is that of the dyadic odometer. As above, let $p_i = (1-x^{2^{i-1}})/2$, so that $s(\brcs{\frac12 (1- x^{2^j})},(\lambda_0 (M_j)) = 0$. On the other hand, a simple computation reveals that the $(1,1)$ entry of $\prod_0^d M_j$ is $2^{-d}\sum_{\eta(a) \equiv 0 \mod n} x^a$ where $a$ varies over the integers from $1 $ to $2^{d+1}$, and $\eta(a)$ is the number of $1$s in the binary expansion of $a$. It easily follows that $\|((1-x)/2)\prod_0^d M_j\|$ is at least $1-1/n - 2^{-d}$ (it can be computed exactly), so that $s_0$ of the sequence is at least $1-1/n$. However, $M_j = M_0(x^{2^j})$, so $s_0 = s_j$ (since the sequence $\brcs{p_i}$ is invariant under $x \mapsto x^2$), so the $s$ value of $(M_j)$ is at least $1- 1/n$. In particular, $(M_j)$ is not isomorphic to the dyadic (or any) odometer.

In fact, we obtain different $s$ values for different values of $n$ (this requires more precise, but elementary computations), so distinct choices of $n$ yield nonisomorphic systems.

\comment
Here is a class of examples yielding nonisomorphisms. Write as previously, the $H$-hemicirculant matrix over $A^+$, $M = \sum q_g m_g$ with $\sum q_g(1) = 1$, and set
$M_j = M(x^{k^j})$ where $ k \geq |H|$, and assume $C:= M(1)$ (evaluating all the entries at $1$) is primitive. Then $(M_j)$ is ergodic (since the underlying real matrix is constant and primitive, hence the real dimension group has unique trace). For a polynomial $q$, $\Log q = \Set{i \in \Z^+}{(q,x^i) \neq 0}$. Now impose the following additional constraints:

\item{(a)} $\brcs{\Log q_g}$ are pairwise disjoint (some may be empty) subsets of the interval $\brcs{0,1,2,\dots, k-1}$;
\item{(b)} setting $d:= \max \Log q_g - \min \Log q_g$, we have $d > 0$ and
$\max \brcs{\cup \Log q} + d < k$.

Now set $p_i = (1-x^{d\cdot k^{i-1}})$ for $i = 1, 2, \dots$. From (b), we deduce (all these are \wrt the usual total variation norm on coefficients) for all $i\in \N$
and $\alpha\in \hat H$,
$$\eqalign{
\left\| p_i \lambda_{\alpha} \(\prod_{j=i}^{j= i+K}M_j \)\right\|& = \left\| p_i \lambda_{\alpha} \(M_i \)\right\|\cr
\left\| p_i \lambda_{0} \(\prod_{j=i+1}^{j= i+K}M_j \)\right\|& = 1 \cr
}$$
Moreover, it is elementary that if $l> i$, then $\| p_l \lambda_{\alpha} \(\prod_{j=i}^{j= i+K}M_j \)\| \geq \| p_i\lambda_{\alpha} \(\prod_{j=i}^{j= i+K}M_j \)\|$. It follows that for each $i$,
$s_i (\brcs{p_t}, (\lambda_{\alpha}(M_j))) = \left\| p_i \lambda_{\alpha} \(M_i \)\right\|$, and because everything in sight is invariant under $x^k$, the latter is simply $ \left\| p_1 \lambda_{\alpha} \(M \)\right\|$, and thus this is the value of $s (\brcs{p_t}, (\lambda_{\alpha}(M_j)))$.

The large eigenvalue, $\lambda_0 (M_j) = \sum q_g (x^{k^{j-1}})$ belongs to $A^+$ and forms the bottom line in ($\dag$). Its $s$ value is thus $ \left\| \frac12 (1-x^d)\sum_g q_g\right\|$.
\endcomment

\comment
Now write $\prod M_j = \sum_{\alpha \in \hat H} v_{\alpha}w_{\alpha} \lambda_{\alpha} (\prod M_j)$ (the range of values in the product consists of a finite interval of positive integers, which need not be specified). Thus $p_i \prod M_j = \sum v_{\alpha}w_{\alpha} p_i\lambda_{\alpha} (\prod M_j)$. It easily follows that $s(\brcs{p_i}, (M_j)) \geq \max_{\hat H} s(\brcs{p_i}, \lambda_{\alpha}(M_j)$, and thus $s(\brcs{p_i}, (M_j)) \geq \max_{\hat H}
\left\| p_1 \lambda_{\alpha} \(M \)\right\|$, and this is $\max_{\hat H}\left\| \frac12 (1-x^d)\sum_g q_g \alpha(g^{-1})\right\|$. Hence, in order to guarantee that $s(\brcs{p_i}, (M_j)) > s\(\brcs{p_i}, \lambda_{0}(M_j)\)$ (and so $(M_j)$ is not isomorphic to the obviously AT system $(\lambda_0(M_j))$---of course, we cannot say whether $(M_j)$ is itself AT), it suffices to show there exists $\alpha$ in $\hat H$ \st $\left\| (1-x^d)\sum_g q_g \alpha(g^{-1})\right\| > \left\| (1-x^d)\sum_g q_g \right\|$ (we have dropped the $\frac12$ factors).

We first note from the definition of $d$ that some cancellation must take place, that is, $ \left\| (1-x^d)\sum_g q_g\right\| < 2$. Next, if $f = \sum c_i x^i$ is a polynomial with nonnegative coefficients and $z_i$ are complex numbers with modulus one, we see
$$\eqalign{
R:=\| (1-x^d)f \| & = \sum |c_i - c_{i-d}|\cr
S:=\| (1-x^d)\sum c_i z_i x^i \| & = \sum |c_i - \frac{z_{i-d}}{z_i}c_{i-d}|\cr
}$$
Since $|a + zb| \geq |a+b|$ with equality only when $z = 1$ (for $a,b > 0$ and $|z| = 1$), we deduce that $S \geq R$ with equality if and only if $c_{i-d} \neq 0$ implies $z_{i-d} = z_i$. However, the definition of $d$ (arising from different choices of elements of $H$) guarantees the strict inequality, $\left\| (1-x^d)\sum_g q_g \alpha(g^{-1})\right\| > \left\| (1-x^d)\sum_g q_g \right\|$. Hence $s(\brcs{p_i},\cdot)$ distinguishes the two systems $(M_j)$ and $(\lambda_0 (M_j))$. \qed

The idea underlying this computation, that we can multiply the diagonal entries (that is, $\lambda_{\alpha} (\prod M_j)$) by $p_i$ to bound $s(\brcs{p_i}, (M_j))$, together with the $S \geq R$ inequality (at least when the supports of $q_g$ are disjoint) obviously appears to be capable of extension to more general situations, possibly even to showing that if $(M_j)$ is an ergodic hemicirculant sequence that is not hollow, then $(M_j)$ is not isomorphic to $(\lambda_0 (M_j))$. What we would really like are conditions on the former guaranteeing that it is not AT or even AT($n-1$), but the mass cancellation invariants are too coarse for this.

In any event, we have shown the following.

\Lem Proposition. Let $M :=\sum q_g m_g$ be an $H$-hemicirculant matrix over $A^+$ \st $\sum q_g(1) = 1$ and $M(1)$ is a primitive (real) matrix. Suppose $M$ satisfies: {\par}
\item{(a)} $\brcs{\Log q_g}$ are pairwise disjoint (some may be empty) subsets of the interval $\brcs{0,1,2,\dots, k-1}$;{\par}
\item{(b)} setting $d:= \min \Set{ (\Log q_g - \Log q_h) \vee 0 }{g \neq h; }$, we have $d > 0$ and
$\max \brcs{\cup \Log q_g} + d < k$.{\par\vskip2pt}
\noindent Set $M_j = M(x^{k^j})$ where $ k \geq |H|$. Then $(M_j)$ is ergodic and {\it not\/} isomorphic to $(\lambda_0 (M_j))$.

The following is practically tautological. Better would be the same conclusion merely hypothesizing $(M_j)$ is not hollow.

\Lem Proposition . Suppose that $(M_j)$ is an ergodic hemicirculant sequence \st $(M_j^2)$ is not hollow. Then $(M_j) \not \iso (\lambda_0 (M_j))$.

\Pf We employ the mass cancellation invariant, using $p_i \in \C[x,x^{-1}]$. Since $(M_j^2)$ is not hollow, there exists $\alpha \in H \setminus \brcs{\chi_0}$ \st for some $N_0$, $\lim_{l \to \infty} \left\|\prod_{j=N_0}^{N_0 + l} \lambda_{\alpha}\( M_j\)^2\right\| > 0$. It follows that given $\epsilon > 0$, there exists $N$ \st $\lim_{l \to \infty}\left\| \prod_{j=N}^{N + l} \lambda_{\alpha}\( M_j\)^2\right\| >1 - \epsilon$.

Set $p_i' = \prod_{j=N}^{N +i} \lambda_{\alpha}\( M_j\)$; then $ 1 \geq \| p_i'\| > 1 - \epsilon$, so we can set $p_i = p_i'/\| p_i'\|$. For each $i$ and all $k$, $ \| \lambda_{\alpha}(\prod_{j=N}^{N+k} M_j)\| > 1-\epsilon$. Hence $s(\brcs{p_i},(M_j)) \geq 1- \epsilon$ for all $\epsilon >0$, and thus $s(\brcs{p_i},(M_j))=1$.

On the other hand, $s({p_i}, (\lambda_0(M_j))) = 0$ from xxx. \qed

\endcomment

Here is a nonisomorphism result, which is almost what we expect: that $(M_j)$ ergodic, hemicirculant, and not hollow implies $(M_j ) \not \iso (\lambda_0 (M_j))$. No counterexamples are known to this conjecture.

\Lem Lemma \foutwo. Suppose that $(M_j)$ is an ergodic sequence of hemicirculant matrices \st $(M_j^2)$ is not hollow. Then $(M_j) \not\iso (\lambda_0(M_j))$. More generally, if $(M_j^{k+1})$ is not hollow, then $(M_j^k) \not \iso (\lambda_0(M_j^k))$.

\Pf If $k=1$, $(M_j^2)$ is not hollow and there thus exists $\alpha \in \hat H\setminus \brcs{\chi_0}$ \st for some positive integer $N_0$, $\lim_{l\to \infty} \left\| \prod_{j= N_0}^{N_0 +l} \lambda_{\alpha}(M_j)^2 \right\| > 0$. It follows that given $\epsilon > 0$, there exists $N \equiv N(\epsilon)$ \st
$\lim_{l\to \infty} \left\| \prod_{j= N}^{N +l} \lambda_{\alpha}(M_j)^2 \right\| > 1-\epsilon$.

Set $p_i' = \prod_{j= N}^{N +i} \lambda_{\alpha}(M_j)$ (in $A_{\C} = \C[x,x^{-1}]$), so that $1 \geq \| p_i'\| > 1 -\epsilon$, and define the normalized versions, $p_i = p_i'/\|p_i' \| \in A_{\C}$. Then for each $i$ and all $k$,
$\left\| p_i \prod_{j= N}^{N +k} \lambda_{\alpha}(M_j) \right\| > 1-\epsilon$. Thus $s\(\brcs{p_i},(M_j)\) \geq s\(\brcs{p_i},(\lambda_{\alpha}(M_j))\) \geq 1-\epsilon$, so that $s\(\brcs{p_i},(M_j)\) \geq 1-\epsilon$ (because the choice of $\brcs{p_i}$ depended on $\epsilon$, we cannot obtain $1$ exactly).

One the other hand, for $N' \geq N$,
$$
p_i \prod_{j= N'}^{N' +k} \lambda_{\alpha}(M_j) = \prod_{i= N}^{N +i} \lambda_{0}(M_j) \prod_{j= N'}^{N' +k} \lambda_{\alpha}(M_j).
$$
If we choose $i$ large enough that $N+ i > N' + k$, then Lemma \onetwo{} applies, so the norm of the product can be made arbitrarily small. Thus $s\(\brcs{p_i},(\lambda_0(M_j)\) = 0$, and therefore nonisomorphism occurs.

If $k > 1$, we find $N$ \st for all $i$, $\left\|\prod_{j=N}^{N+i}\lambda_{\alpha}(M_j)^{k+1}\right\| > 1-\epsilon$, set $p_i = \prod \lambda_{\alpha} (M_j)$, and proceed as above.\qed

There do exist non-hollow $(M_j)$ with $(M_j^2)$ hollow, \st the same nonisomorphism holds, so the question remains whether the additional hypothesis really is necessary.

Suppose $M_j = \frac 12 \(\smallmatrix 1 & x^{2^j} \\ x^{2^j} & 1 \\ \endsmallmatrix\)$. Radu Monteanu has shown (unpublished) that this is represents the Morse-Thue transformation. It was shown by Dooley \& Quas [DQ] to be AT. It is easy to verify that $M_j^2$ is hollow: $\left\|\prod_{j = N}^{N+k-1} (1-x^{2^j})^2 \right\| = \oh{4^{k}}$ (in fact, it is $\oh{\beta^k}$ for $\beta$ substantially less than $4$). However, $(M_j) \not\iso \( \frac12 (1+x^{2^j})\) = (\lambda_0(M_j))$ as the Morse is not an odometer; we can use $p_i = \frac12 (1 - x^{2^j})$).

The map $AH \to AH$ given by $e_g \mapsto e_{g^{-1}}$ induces an isomorphism between the dimension spaces of $(M_j)$ and of $(M_j^T)$. This allows us to show that not only is the map in Lemma \throne{} not generally an isomorphism between $(M_j \otimes M_j^T)$ and $(M_j M_j^T)$ (which we already knew), but in fact, the two dimension spaces are not generally isomorphic.

\Lem Corollary \fouthr. Suppose that $|H|$ is odd and $(M_j)$ is an ergodic hemicirculant sequence \st $(M_j^2)$ is not hollow. Then $(M_j \otimes M_j^T) \not \iso (M_j M_j^T)$.

\Pf First, $(M_j^T) \iso (M_j)$. Next, $(M_jM_j^T)$ is symmetric and of odd size, so is hollow (Proposition \twothr). If $(M_j \otimes M_j^T) \iso (M_j M_j^T)$, then we would have the following chain of isomorphisms,
$$
(M_j^2) \iso (M_j \otimes M_j) \iso (M_j \otimes M_j^T) \iso (M_j M_j^T) \iso (\lambda_0(M_jM_j^T)) = (\lambda_0 (M_j^2)).
$$
This contradicts $(M_j^2) \not \iso (\lambda_0 (M_j^2))$.
\qed

\Lem Example. Another computation.

\noindent Let $ H = \Z_n$ be the cyclic group of order $n$, let $P$ be the corresponding cyclic permutation matrix of size $n$ (corresponding to the generator, $[1]$ of $\Z_n$, and for $j =0, 1,2, \dots$, set $M_j = \frac 12 (\I + x^{3^j}P)$. Then $(M_j)$ and $(M_j^2)$ are ergodic but not hollow, but $(M_j^3)$ is hollow.

It is reasonable to conjecture that $(M_j)$ is not AT, or even not AT($2$); but by Theorem \onethr(a), $(M_j^2)$ is AT. We wish to determine the corresponding sequence of polynomials, $(Q_j)$ (each $Q_j \in A^+$ and $Q_j(1) = 1$) that is obtained by the argument in op cit \st $(M_j^2) \iso (Q_j)$.

Indexing $\hat H $ by $0, 1, 2, \dots, n-1$ (instead of $\alpha$), we see immediately that $\lambda_k (M_j) = \frac12(1 + x^{3^j}\xi^k)$, where $\xi = \exp (2\pi i/k)$, for $k = 0,1, \dots, n-1$. For $k \neq l$ (corresponding to $\alpha \neq \beta$), we calculate the norm of $\lambda_k (M_j) \lambda_l (M_j)$. In all cases, we see that the maximum of the norm over all pairs is some number $\gamma \equiv \gamma (n) < 1$. As a sample computation, $\gamma(3) = \frac{3}{4}$, comes from
$$\eqalign{
(1+ \xi x)(1 + \xi^2 x) &= 1 - x + x^2\cr
(1+ \xi x)(1 + x) &= 1 - \xi^2 x + \xi x^2\cr
}$$
(because $M_j = M_0 (x^{3^j})$, the dimension space is stable under $x \mapsto x^3$ and we need only consider the case that $ j= 0$; also, when $n=3$, $1 + \xi + \xi^2 = 0$).

Now we find a telescoping $0 = n(1) < n(2) < n(3) <\dots $ so that for each $k \neq l$,
$$\sum_{t}\left\|\lambda_k \(\prod_{n(t)}^{n(t+1)-1} M_j\)\lambda_l \(\prod_{n(t)}^{n(t+1)-1} M_j\)\right\| < \infty;$$ sufficient for this is simply
$$
\sum \gamma^{n(i+1) -n(i)} < \infty.
$$
We can set $n(i+1) = i(i+1)/2$ (so that $n(i+1) -n(i) = i$---this is a lot bigger than necessary. Now define $M^{(i)} = \prod_{n(i)}^{n(i+1)-1} M_j$. We find the coefficients in the expansion of $M^{(i)} = \sum_{u \in H} c_u m_u$ (as a circulant matrix).

For an integer $t$, let $t = \sum \eta_i(t)3^i$ be its ternary representation ($\eta_i(t) \in \brcs{0,1,2}$); let $V_m = \Set{t < 3^{m+1}}{\eta_i(t) \in \brcs{0,1} \text{ for all $i$}}$; and for $t $ in $V_m$, define $\eta(t) = \sum_i \eta_i (t)$, that is, the number of $1$s in the ternary expansion of $t$. For $ j \in \hat H$ and $ 0 < a < b$, define
$$
P_{j, a, b}:= \sum_{\Set{t \in V_{b-a}}{\eta (t) \equiv j \mod n}} x^{t3^a}.
$$
Then a simple induction argument yields that with $c_{ui} = P_{u,n(i),n(i+1)}$, we have $M^{(i)} = \sum_{u \in H} c_{ui} P^u/2^i$, using $n(i+1) - n(i) = i$ (this can also be obtained by noting that $\lambda_u (M^{(i)}) = \prod (1 + \xi^u x^{u 3^{n(i)}})/2^i$ and applying the finite Fourier transform).

By the method of Theorem \onethr, we have an approximate factorization $(M^{(i)})^2 \sim V^{(i)} W^{(i)}$, where $V^{(i)} = \frac1n\sum e_u c_{ui}$ and $W^{(i)} = \sum e_u c_{n-u,i}$. This yields an isomorphism between $(M_j^2)$ and $(p_i: = W^{(i+1)}V^{(i)})$. The earlier computation yields
$p_i = \sum_{u \in \Z_n} P_{u,n(i),n(i+1)}P_{n-u,n(i+1),n(i+2))})$, and this equals $P_{0,n(i),n(i+2)}$. Thus $(p_i = P_{0,n(i),n(i+2)})$. It is very likely that $((\frac 1{|H|} \tr (M^{(i)}))^2)$ is not isomorphic to $(p_i)$, but at the moment, there are technical difficulties in showing this. (The idea is to use the mass cancellation invariants discussed above.)

In this example, all of the automorphisms on $(M_j^2)$ given by $H$ are nontrivial. This is one way of constructing finite order automorphisms of AT actions that induce nontrivial automorphisms of the dimension space. Of course, this realizes the dual action of the original approximately inner action of $H$ on the underlying ITPFI given by $(\lambda_0(M_j^2))$.\qed

\Lem Example. An example of a hemicirculant ergodic sequence $(M_j)$ \st $(M_j \otimes M_j^T)$ is not isomorphic to $(M_j M_j^T)$.

\noindent This is a modification of the previous example. Let $H = \Z_n$ where $n \geq 5$ and is odd. Select $k \geq 4$, and set $M_j = \frac12\(\I + Px^{k^{j}}\)$. Since $k\geq 4$, $(M_j^3)$ is not hollow, and the result follows from Corollary \fouthr.
\qed

If $(M_j)$ is hollow, there exists a telescoping $n(i)< n(i+1)$ so that if $M^{(i)}:= M_{n(i+1)-1}\cdot M_{n(i+1)-2}\cdot \dots \cdot M_{n(i)}$, then $(M_j) \iso (\tr (M^{(i)}))$. However, we noted that for some sequences of hemicirculant matrices of the form $(M_j^2)$, the isomorphism is only obtained with overlapping sequences;
that is $(M_j^2) \iso (\tr (M^{(i+1)}M^{(i)}))$ (and generally not isomorphic to $\(\tr ((M^{(i)})^2)\)$). This overlapping is essential (and one might view it as natural, since $(M_jM_{j-1}) \iso (M_j^2)$), and extends to the most general case.

Suppose $(\Arrow N_i; A^{m(i)}.A^{m(i+1)})$ is a sequence of matrices with entries from $A^+$ \st the column sums, when evaluated at $x=1$, are all $1$. Suppose it is known that $(N_j)$ (or more accurately, its associated dimension space) is AT. Then, after a possible telescoping (which we incorporate into the notation, to avoid unnnecessary complicating an already complicated situation), there exist columns $V_i$ (in $(A^{m(i+1) \times 1})^+$) and columns $W_i$ (in $(A^{1 \times m(i) })^+$)
\st $\sum \tripnorm N_i - V_i W_i xxx < \infty$. Then $(N_i ) \iso (W_{i+1}V_i)$, the latter of course being a sequence of normalized elements of $A^+$ (we can slightly perturb the sequence so that the values at $x = 1$ are all exactly one, rather than close).

Now take any telescoping of the (already telescoped) $(N_j)$, say corresponding to $0 = n(1) < n(2) < \dots$, and this time, take the products including {\it both\/} ends of the intervals; that is, define $N_{(i)} = N_{n(i+1)}\cdot N_{n(i+1)-1}\cdot \dots \cdot N_{n(i)}$ (the more usual telescoping would not include $N_{n(i+1)}$). Set $P_j = W_{j+1}V_j$ (an element of $A^+$). Then
$$\eqalign{
N_{(i)} &\sim V_{n(i+1)}W_{n(i+1)}V_{n(i+1)-1} W_{n(i+1)-1} \dots V_{n(i)}W_{n(i)} \cr
& = V_{n(i+1)} W_{n(i)} \( W_{n(i+1)}V_{n(i+1)-1} \cdot\dots \cdot W_{n(i)+1}V_{n(i)} \) \cr
& = (V_{n(i+1)} W_{n(i)}) \cdot P_{n(i+1)-1} \cdot P_{n(i+1)-2} \cdot \dots \cdot P_{n(i)} = (V_{n(i+1)} W_{n(i)}) \prod_{j=n(i)}^{n(i+1)-1} P_j ;\cr &\quad \text{thus, if $N_j$ were square matrices,} \cr
\tr N_{(i)} & \sim (W_{n(i)}V_{n(i+1)})\prod_{j=n(i)}^{n(i+1)-1} P_j.\cr
}$$
Since $(N_j) \iso (P_j)$, we obtain the isomorphism of AT systems, $(\tr N_{(i)}) \iso (N_j)\otimes (W_{n(i)}V_{n(i+1)}))$. This is true for {\it every\/} telescoping of the originally telescoped sequence $(N_j)
$ (although in principal, the isomorphism class of $(\tr N_{(i)} )$ may depend on the choice of telescoping, $n(i)$). We note in particular, the overlapping phenomenon---$N_{(i)}$ and $N_{(i-1)}$ both contain the term $N_{n(i)}$ (in the former, at the extreme right of the product, in the latter, at the extreme left), although no such overlapping occurs with the $P_j$s.

One might expect that if the telescoping were sufficiently sparse (that is, if $i\mapsto n(i) $ grows really quickly), the extra polynomials $W_{n(i)}V_{n(i+1)}$ should not play much of a role, that is, it is possible that $(\tr N_{(i)}) \iso (N_j)$ for sufficiently fast growing $n(i)$.

The overlapping has significance for actual computations; for example, when $N_j$ is hemicirculant, it is often straightforward to calculate the $(1,1)$ entry of $N_{(i)}$ (this is the trace divided by $n$). However, it is not so straightforward to calculate even $\tr N_{(i)}N_{(i-1)}$ (and computing the traces of products of larger numbers of matrices is still more tedious), precisely because of the overlap.

It also means that a theorem such as, if $(N_j)$ is AT, then there is a telescoping so that $\tr (N_{n(i)} N_{n(i-1)})$ is almost $\tr N_{n(i)} \tr N_{n(i-1)}$ (i.e., the trace is almost multiplicative after suitably telescoping), is almost certainly not true. This would have been really useful to show that many sequences are not AT.

\SecT 5 Perspective

\noindent In terms of the type III von Neumann algebras that are classified by dimension spaces, families of hemicirculant matrices have a particularly elementary interpretation. Begin with a product type W*-algebra (and corresponding AT action and dimension space), let $H$ be a (finite abelian) group of product type automorphisms (necessarily approximately inner); assume the action is ergodic, and form the crossed product, another W*-algebra. Its dimension space is given by an ergodic sequence of $H$-hemicirculant matrices, and the inclusion of the original in the crossed product translates in the dimension space setting to the diagram (*) (or whatever it is).

The original action of $H$, being approximately inner, is not otherwise visible on the underlying AT dimension space; however, its dual action on the crossed product is represented by the action of $H$ in its regular representation as permutation matrices; these commute the $(M_j)$. Of course, the dual action may also turn out to be approximately inner, in which case these permutation matrices act trivially on the dimension space level; this occurs precisely when the sequence $(M_j)$ is hollow, as is easy to verify.

In terms of $\Z$-actions (i.e., ergodic transformation), sequences of hemicirculant matrices correspond to a class of $H$-actions on an AT dynamical system (presumably, there are non-product type actions inequivalent to any product action, even for $H = \Z_2$, hence not all $H$-actions are determined by hemicirculant matrices), resulting in an ae finite to one (at most $|H|$ to one) map between the spaces. (A particular consequence is that the resulting systems have entropy zero, but it is likely that AT($n$) implies entropy zero anyway.) Hollow sequences yield one to one maps, that is, isomorphisms.

We can also relate this to corresponding ideas from topological (rather than measure-theoretic) classification, in this case of actions of compact groups on AF C*-algebras (in particular). From an AF algebra, given as a direct limit of finite-dimensional algebras, $ \lim A_k$ dense in $C$, let $G$ be any compact group (not necessarily abelian), and suppose for each $k$, there is a unitary representation $\Arrow \pi_k ; G.A_k $ compatible with the maps $A_k \to A_{k+1}$. This yields an action of $G$ on $A$ (called {\it locally representable\/}).

It was shown in [HR, Theorem III.1], that $K_0 (C \times G) $ (viewed as an ordered module over the representation ring) together with some additional data is a complete invariant (up to unitary group actions on C*-algebras) for the action. Moreover, all such actions can be obtained by considering all possible Bratteli diagrams that represent $C$, and for each one, replacing the integer entry (describing the multiplicity) by a character (not necessarily irreducible) of $G$ whose dimension equals the multiplicity (that is, if $\chi$ is a character, the arrow with multiplicity $d$ can be weighted by $\chi$ if $\chi(1) = d$. This yields Bratteli diagrams with weights from the positive cone of the representation ring.

If we begin with a Bratteli diagram that is initially weighted with polynomials (in one variable) with positive integer coefficients (necessarily not adding to one), we obtain an action of the circle on the underlying AF algebra. For example, a product type action would correspond to a sequence of polynomials in one variable. We can further elaborate on this diagram by now permitting characters of $T \times H$ (where $H$ is a finite abelian group), yielding an $H$-action on the crossed product by $T$. We can then convert this into a product type action in the measure theoretic setting selecting an extreme trace on the underlying AF algebra and then dividing by suitable rationals so that things add to one.

Going in reverse, we merely have to approximate sufficiently well any real numbers that appear in the coefficients in the matrix entries by rational numbers (since the measure-theoretic classification is immune to tiny perturbations) and then multiplying each matrix entry by a suitable positive integer to remove the denominators. Of course, measure-theoretic classification is much coarser than topological, and topological classification is sensitive to the approximation and choice of integer.

\vskip 10 pt
\SecT 6 Powers of transformations

Let $(T,X,\mu)$ be a dynamical system, with $T$ ergodic and $\mu$ quasi-invariant \wrt $T$. The Poisson boundary can be represented via the dimension space of a suitable sequence of maps $\Arrow N_j; A^{m(j)}. A^{m(j+1)}$ together with a trace induced factoring through an extremal trace on the direct limit of $\Arrow N_j(1); \R^{m(j)}. \R^{m(j+1)} $ where evaluation at $1$ means evaluation $x \mapsto 1$ to each coordinate, creating a real column stochastic matrix.

Suppose we want to consider the square or a higher power of $T$, that is, the dynamical system $(T^n,X,\mu)$ for some $n > 0$. A first problem is that this need not be ergodic (the use of dimension spaces is more or less confined to ergodic transformationsthese are insensitive to sufficiently small perturbations, whereas in the non-ergodic case, arbitrarily small perturbations can add or delete atoms). However, in general there is a process for obtaining a dimension space for $T^n$ (where convenient, we abbreviate the triple $(S,X,\mu)$ to its first component) out of the original one for $T$.

We can read off from the construction whether $T^n$ is ergodic, and if so, the resulting dimension space corresponds to a special $\Z_n$ action; in particular, if $T$ is AT, then the dimension space for $T^n$ comes from a commuting sequence of $n \times n$ matrices that are very close to being circulant, and allow us to read off properties. For example, if $T$ is AT($k$) and $T^n$ is ergodic, then $T^n$ is AT$(kn)$ (but it frequently happens that it is AT). We give examples wherein $T$ is AT, $T^2$ is ergodic, but $T^2$ is not AT. It also follows from earlier results that if $T^n$ is ergodic and $T$ is AT, then $T^n \otimes T^n$ is AT.

First we give the construction. Let $X = x^n$, and form $B = \R [X,X^{-1}]$ sitting inside $A = \R[x^{\pm1}]$, each equipped with the usual ordering. Then $A$, viewed as a $B$-module, is free on the set $\brcs{1,x,\dots,x^{n-1}}$. Write $A = \oplus A_i$ (as $B$-modules) where $A_i = x^i B$. Now let $p$ be a polynomial in $A$ with no negative coefficients. Let $Q$ be the companion matrix of the polynomial (in $Z$), $Z^n - x^n = Z^n - X$,
$$
Q = \(\matrix 0 & 0 & 0 & \dots & 0 & X\\
1 & 0 & 0 & \dots & 0 & 0\\
0 & 1& 0 & \dots & 0 & 0\\
&& \ddots &&&\\
0 & 0 & 0 & \dots & 1 & 0\\
\endmatrix \).
$$
Then $Q$ belongs to M$_n B^+$; form $\Cal B (p):= p (Q)$. Then $p(Q)$ has entries in $B^+$, its column sums when evaluated at $x =1$ (and at $X =1$!) equal $p(1) = 1$, and if we form the diagonal matrix (with entries from $A^+$), $\Delta = \diag (x^{n-1},x^{n-2}, \dots, 1)$, we have $\Delta Q \Delta^{-1} = xP$, where $P$ is the cyclic permutation matrix. In particular, $\Delta \Cal B(p) \Delta^{-1} = p(xP)$, a circulant matrix. Of course, this is not implementable over $B$, so we cannot simply transfer everything in what follows to circulant matrices.

Now the map $\Arrow \times x;A.A$where we identify $A$ with the free $B$ module, $A = \oplus_{i=0}^{n-1} x^i B$has $Q$ as its matrix representation (\wrt the basis $\brcs{x^i}$). There is a natural map $\Arrow \varphi;B^n. A$ sending $(b_0, b_2, \dots,b_{n-1})^T \mapsto \sum x^i b_i$ (which effectively is the identity on $A$ after identifications).

Now let $N_j$ be an $m(j+1) \times m(j)$ matrix with entries from $A^+$ \st all the column sums of $N_j (1)$ (evaluation at $x=1$) are $1$. Form a matrix $\Cal B(N_j)$ of size $m(j+1)n \times m(j)n$ by replacing each entry $(N_j)_{lm}$ of $N_j$ by the $n \times n$ matrix $\Cal B((N_j)_{lm}) = (N_j)_{lm} (Q)$. The resulting matrix has entries in $B^+$, and its column sums, when evaluated at $X=1$ are all $1$. The claim is now that the dimension space (now viewed as a $B$-module, equivalently as a module over $l^1 (n\Z)$) corresponding to $(\Cal B(N_j))$ implements $T^n$.

This is practically a tautology. The natural map $\Arrow \varphi;B^n.A$ extends to $(B^n)^m \to A^m$ in the obvious way, and the diagram
$$
\diagram
(B^{n})^{m(j)}&\rTo^{{\Cal B}(N_j)\hphantom{xx}} & (B^{n})^{m(j+1)}&\rTo^{{\Cal B}(N_{j+1})} &
(B^{n})^{m(j+2)}&\rTo^{{\Cal B}(N_{j+2})} & \dots & \\
\dTo^{\varphi} &&\dTo^{\varphi}& & \dTo^{\varphi} & && \\
A^{m(j)}&\rTo^{N_j} & A^{m(j+1)}&\rTo^{ N_{j+1}} &
A^{m(j+2)}&\rTo^{N_{j+2}} & \dots & \\
\enddiagram
$$
commutes and the vertical maps are isomorphisms of ordered $B$-modules. Under this, multiplication by $X$ in the terms of the top row translates to multiplication by $x^n$ on the bottom row. Hence $T^n$ is represented by the $B$-dimension space arising from $(\Cal B(N_j))$. An obvious consequence is that if $T$ is AT($k$) and $T^n$ is ergodic, then $T^n$ is AT($nk$) (the definition of AT$(m)$ does not require ergodicity; with ergodicity, AT$(m)$ is equivalent to the sequence of matrices being at most $ m\times m$, whereas without ergodicity, only the reverse implication applies).

At this stage, we may be tempted to use $\Delta$ to change this to matrices blown up by circulant matrices. If we can choose $m(j)=1$ for all $j$, that is, if $T$ is AT, after applying $\Delta$, we would obtain a sequence of circulant matrices, and so use all the results of the preceding section. Unfortunately, this is not permissible, since we are restricted to $B$- (or $l^1 (n\Z)$-) maps between dimension spaces. Fortunately, we can often deduce results from the circulant case.

The test for ergodicity, however, is identical, since after evaluation at $X=1$, the resulting $\Cal B(p_j)(1)$ is circulant. For example, if $(p_j = (1+x^{2^{j-1}})/2)$, then $T$ is the dyadic odometer. Obviously $T^2$ and thus $T^{2n}$ are not ergodic$T^2$ is the disjoint union of two copies of the odometerwe can also see this from $\Cal B (p_j) = \frac 12 \diag (1+ X^{2^{j-1}},1+ X^{2^{j-1}}) $, so the limiting dimension group is just a direct sum of two copies the original odometer sequence (with $X$ playing the role of $x$).

On the other hand, if $n$ is odd, $T^n$ is ergodic; $\Cal B(p_j)$ is $(\I + Q^{2^{j-1}})/2$, which when evaluated at $X=1$ yields the sequence of matrices $\frac12 (\I + P^{2^{j-1}})$. Since $n$ is odd, there exist infinitely many $k$ \st $2^k \equiv 1 \mod n$, and it easily follows that the matrix products converge to the standard rank one matrix. In fact, for odometers, there is a much stronger result available. As Thierry Giordano observed, any ergodic power of an odometer is conjugate to the original, as follows from the ($L^1$) pure point spectrum property; we provide a much more complicated proof of the same result to illustrate the methods.

It can also be derived directly from the dynamical origins, at least in the case that the supernatural number is of the form $n^{\infty}$ for some $n$, and likely more generally. Let $Y = \prod_{\N} \Z_n$ with the product topology (a Cantor set); the odometer views $Y$ as $\Z_{(n)}$ (the $n$-adic completion of the integers; the elements of the sequence space can be viewed as power series), and simply adds $1$. Assume $(k,n) =1$; the $k$th power of the odometer adds $k$. However, since $k$ is relatively prime to $n$, $k$ is invertible in $\Z_{(n)}$ and thus $\times k$ (multiplication by $k$) is also a self-homeomorphism of $Y$, using the multiplicative (rather than the additive) structure of $\Z_{(n)}$. Then we have the commuting diagram
$$
\diagram
Y&\rTo^{+1}& Y\\
\dTo^{\times k} & &\dTo^{\times k} \\
Y&\rTo^{+k}& Y\\
\enddiagram
$$
The top row is the odometer, the bottom row is its $k$th power, and the vertical maps are invertible maps (exploiting the multiplicative structure) which conjugate the odometer to its $k$th power. (It is not clear how to extend this elementary argument to more general odometers for which there is no underlying multiplication, e.g., if the supernatural number has no infinite multiplicities.)

\comment

If $k$ is a positive integer, then we say {\it $k$ is prime to a supernatural number $\Arrow h; \text{Spec \Z}. \N \cup \brcs{\infty}$\/} if $h(p) =0$ for all primes $p$ dividing $k$.

\Lem Proposition \sixone. If $S$ is an odometer and $k$ is relatively prime to the supernatural number of $S$, then $S^k$ is conjugate to $S$.

\Pf Let $X = x^k$, $B = \R[X^{\pm1}]$, and form $Q$ as before, the companion matrix of $\lambda^k - X$. Let $(n(j))$ determine the odometer (in the sense that the prime divisors in $n(j)$ occur with total multiplicities over all $j$ agreeing with the supernatural number), set $T(j) = \prod_{i=1}^{j} n(i)$, and $n(1) = 1$. All $n(j) $ are relatively prime to $k$. By telescoping, we may assume that $T(j) \equiv 1 \mod k$ (and thus $n(j) \equiv 1 \mod k$), and $\sum 1/n(j) < \infty$. Form $f_j:= \sum_{s = 0}^{n(j+1)-1} (x^{T(j)})^s$, so that the dimension space of the odometer is determined by $(f_j/n(j+1))$.

Fix $j$ and let $r(i)$ denote the remainder on dividing $i$ by
$k$ (so that $r(i) \in \brcs{0,1,2,\dots,k-1}$). Now we consider the effect of the multiplicative map $x\mapsto Q$ on $f_j$. First, $x^{sT(j)} \mapsto Q^{sT(j)} = X^{\flo{sT(j)/k}}Q^{r(sT(j))}$. Since $T(j) \equiv 1 \mod k$, the power of $Q$ appearing in this expression is simply $r(s)$. Hence the image of $f_j$ under this map is
$$
f_j \mapsto \sum_{s = 0}^{n(j+1)-1} X^{\flo{sT(j)/k}}Q^{r(s)}.
$$
Now for each $a \in \brcs{0,1,\dots,k-1}$, the coefficient of $Q^a$ in this expansion is the sum
$$\sum_{s \equiv a \mod k, 0 \leq s \leq n(j+1)-1} X^{\flo{sT(j)/k}}.
$$
Write $T(j) = kl +1$ (where $ l \equiv l(j) = (T(j)-1)/k$) and for $s$ appearing in the sum, set $s = kt + a$. Then $sT(j) = a + k(t +al + ktl)$, and thus $\flo{sT(j)/k} = t + al + ktl = t(kl+1) + l = tT(j) + al$. As $s$ varies over the set, $t$ varies over the interval $\brcs{0,1,\dots, \flo{(n(j+1)-a-1)/k}}$, and we note that as $a$ varies over $\brcs{0,1,\dots, k-1}$, the upper limit for the value of $t$ changes by only $1$. Set $X_j = X^{T(j)}$ and $P_j = \sum_{t=0}^{\flo{(n(j+1)-1)/k}} X_j^t$, so that the error introduced (after normalizing by dividing by $f_j(1) = 1/n(j+1)$) in approximating the coefficient of $Q^a$ by $X^{al}P_j$ is at most $1/n(j+1)$. Hence
$$
\left\| \Cal B\(\frac{f_j}{n(j+1)} \)- \sum_{a=0}^{k-1} (QX^l)^a \frac{P_j}{n(j+1)} \right\| < \frac k{n(j+1)}.
$$
Now we analyze $\sum_{a=0}^{k-1} (QX^l)^a$. A simple computation reveals that if $\Phi_j \equiv \Phi = \diag(X^{(k-1)l}, X^{(k-2)l}, \dots , X^l , 1)$, then (recalling that $T(j) = kl+1$, and $X_j = X^{T(j)}$)
$$
\Phi \sum_{a=0}^{k-1} (QX^l)^a \Phi^{-1} = \(\matrix 1 & X_j & X_j & \dots &X_j & X_j \\
1 & 1 & X_j & \dots &X_j & X_j\\
& & \ddots & & & \\
1 & 1 & 1 & \dots &1 & X_j \\
1 & 1 & 1 & 1 \dots &1 & 1\\
\endmatrix\).
$$
The entries in the strict upper half are all $X_j$, while those in lower half are all $1$. Now $\| (X_j P_j -P_j)/n(j+1)\| = 1/n(j+1)$ (since multiplication by $X_j$ merely shifts the powers of $X_j$ in $P_j$ by one). Since we are taking the column sum norm, we deduce (with $\pmb 1$ denoting the matrix all of whose entries are $1$)
$$\eqalign{
\left\|\Phi P_j \sum_{a=0}^{k-1} (QX^l)^a \Phi^{-1} - \frac{P_j}{n(j+1)/k} \pmb 1\frac 1k\right\| &< \frac k{n(j+1)}\text{\/; \qquad and thus }\cr
\left\| P_j\sum_{a=0}^{k-1} (QX^l)^a - \frac{P_j}{n(j+1)/k}\Phi^{-1} \pmb 1\frac 1k \Phi\right\| &< \frac k{n(j+1)}\cr
}$$

Set $V_j = (X^{-(k-1)l}, X^{-(k-2)l}, \dots, 1)\frac{P_j}{n(j+1)/k}$ and $W_j = \frac1k (X^{(k-1)l}, X^{(k-2)l}, \dots, 1)^T$. The last inequality then yields
$$\eqalign{
\left\|P_j\sum_{a=0}^{k-1} (QX^l)^a - W_j V_j \right\|&< \frac k{n(j+1)}\text{\/; \qquad and thus } \cr
\left\|\Cal B(f_j) - W_j V_j\right\| &< \frac {2k}{n(j+1)}.
}$$
Since we have assumed that $\sum 1/n(j) < \infty$, this shows that $T^k$ has AT dimension space (over $B$) implemented by the sequence of polynomials $(V_{j+1}W_j)$. Each $V_{j+1}W_j$ expands to $P_{j+1}/((n(j+2)-1)/k)\cdot \sum_{a=0}^{k-1} X^{a(l(j)-l(j+1))}/k$. Denote the second factor by $p_{j+1}$; then
$$\eqalign{
p_{j+1} &= \frac{\sum_{a=0}^{k-1} X^{a(l(j)-l(j+1))}}{k} \cr
&= \frac{\sum_{a=0}^{k-1} (X^{-(n(j+1)-1)T(j)/k})^a}{k} \cr
& = \frac{X^{-(n(j+1)-1)T(j)}\sum_{a=0}^{k-1} (X^{(n(j+1)-1)T(j)/k})^a}k.
}$$
The second line comes from $l(j+1)- l(j) = (n(j+1)T(j) - 1 + T(j)-1)/k$, noting that $n(j+1) \equiv 1 \mod k$. Now
$$\eqalign{
P_j p_{j+1} & = \frac{1}{k}X^{-(n(j+1)-1)T(j)}\sum_{s=0}^{\flo{(n(j+1)-1)/k}} X^{T(j)s} \cdot \sum_{a=0}^{k-1} (X^{(n(j+1)-1)T(j)/k})^a\cr
& = \frac 1k X^{-(n(j+1)-1)T(j)} \sum_{s=0}^{n(j+1)-1} X^{T(j)s}\cr
& = \frac 1k X^{-(n(j+1)-1)T(j)} n(j+1)f_j(X).\cr
}$$

Since we have isomorphisms (of $B$-dimension spaces) $(\Cal B(f_j)) \iso (\frac {P_{j+1}}{(n(j+2)/k} \cdot p_{j+1}) \iso (\frac {P_{j+1}}{(n(j+2)/k})\otimes (p_{j+1})$, we also deduce the isomorphism obtained by shifting the $P_j$, that is, $(\Cal B(f_j)) \iso (\frac {P_{j}}{(n(j+1)/k})\otimes (p_{j+1}) \iso (\frac {P_{j}}{(n(j+1)/k} \cdot p_{f+1}) = (X^{-(n(j+1)-1)T(j)}\cdot f_{j}f_j (X) ) \iso (f_j(X))$ (powers of $X$ can be thrown in randomly without affecting the isomorphism class). The last of course is just the dimension space of the original odometer, this time in $X$.
\qed
\endcomment

\comment
\Lem Proposition. Let $\(\Arrow M_j;B^{d(j)}. B^{d(j+1)}\)$ and $\(\Arrow N_j;B^{e(j)}. B^{e(j+1)}\)$ be ergodic sequences (as $B$-modules). If $\(\Arrow M_j; A^{d(j)}. A^{d(j+1)}\)$ is isomorphic to $\(\Arrow N_j; A^{e(j)}. A^{e(j+1)}\)$ as ordered $A$-modules, then $\(\Arrow M_j;B^{d(j)}. B^{d(j+1)}\) \iso \(\Arrow N_j;B^{e(j)}. B^{e(j+1)}\)$. If $d(j)$ is bounded and $\(\Arrow M_j; A^{d(j)}. A^{d(j+1)}\)$ is AT$(m)$, then $\(\Arrow M_j;B^{d(j)}. B^{d(j+1)}\)$ is AT($m$).

\noindent {\it Remark} In other words, if two $B$- (or $l^1(k\Z)$-) dimension spaces become isomorphic on tensoring over $B$ with $A$ (or $l^1 (\Z)$, then they were isomorphic to begin with.

\Pf fill me in.
\endcomment

\SecT 7 Unexpected isomorphisms

We would expect that the following is either false, or it is true but easy to prove. It is true but (currently) is not easy to prove.

\Lem Proposition \sixthr. Suppose $(M_i)$ is an ergodic sequence of matrices over $B^+$. 
\item{(a)} Suppose that $(M'_i)$ is an ergodic sequence of matrices over $B^+$ such that there exists an $A^+$-module isomorphism between the corresponding dimension modules; then they are isomorphic as $B^+$-modules, that is, $(M_i)_B \iso (M'_i)_B$.
\item{(b)} If there is an $A^+$-module isomorphism $(M_i) \iso (N_i)$ where $(N_i)$ is an AT$(n)$ ergodic sequence with entries from $A$, then $(M_i)$ is AT$(n)$ \wrt $B$.

In particular, If viewed as an ergodic sequence of matrices with entries from $A^+$, $(M_i)_A$ is AT, then $(M_i)_B$ is AT. Moreover, if $(M_i)_A \iso (p_i)_A$ where $p_i \in B^+$ (and $p_i(1) = 1$), then $(M_i)_B \iso (p_i)_B$. This justifies (in some cases) using the conjugation argument available over $A$ but not over $B$ to determine isomorphism over $B$.

\comment
\Lem Proposition \sixthr. Suppose $(M_i)_B$ is an ergodic sequence of matrices over $B^+$ \st every column sum of $M_i (1)$ is $1$. If viewed as an ergodic sequence of matrices with entries from $A^+$, $(M_i)_A$ is AT, then $(M_i)_B$ is AT. Moreover, if $(M_i)_A \iso (p_i)_A$ where $p_i \in B^+$ (and $p_i(1) = 1$), then $(M_i)_B \iso (p_i)_B$.

\endcomment

\noindent {\it Elementary properties of $A = \R[x^{\pm1}]$, $B = \R[x^{\pm k}]$ equipped
with the
$l^1$ norms.}
\def\proj #1{{}^{(#1)}}

\def\dist{\text{dist}\,}

\noindent Write $a = \sum_{i\in \Z_k} x^i a\proj i$ in terms of the standard basis of $B$ as an $A$-module, $\brcs{1,x,\dots, x^{k-1}}$. Then

\item{} $\| a \| = \sum \| a\proj i\|$
\item{} $\dist (a,B) = \sum_{i \in \Z_k \setminus \brcs{0}
} \| a\proj i\|$ and is achieved at $b = a\proj 0$
\item{} if $a \in A^+$, then $\| a \| = a(1)$ and all $a\proj i \in B^+$
\item{} if $a \in A^+$, then $\dist (a,B) = \dist (a,B^+)$
\item{} if $a_j \in A^+$, then $\dist ( \sum a_j,B) = \sum \dist (a_j, B)$


\Lem Lemma \sixtwo. Suppose $a,c$ are in $A^+$ and there exists $b$ in $B^+$ \st
$\| ac - b\| < \epsilon a(1) c(1)$. Provided $ (k^2 - k) \epsilon < 1/2$,
there exists $j $ in $\brcs{0,1,2,\dots, k-1}$ together with $a'$ and $c'$
in $B^+$ \st $\| x^j a' - a\| < (k^2-k)\epsilon a(1)$ and $\| x^{-j}c' - c\|
<(k^2 - k) \epsilon c(1)$.

\Pf Write $a = \sum_{i=0}^{k-1} a_i x^i$ and $c = \sum_{i=0}^{k-1} c_i
x^i$ where $a_i, c_i \in B^+$. For some real number $r > 0$, suppose that
at least two of $a_i$ have norm at least $r a(1)$. Select $l$ \st $c_l (1)
\geq c(1)/k$. Because all the terms have coefficients that are in $B^+$,
the norm is just evaluation at $x = 1$ or $X = 1$. Of the (at least) two
choices for $i$ \st $a_i (1) > ra(1)$, at least one of them, say $i'$, is
not congruent modulo $k$ to $-l$. Hence in the product $ac$, the term $x^{i' +
l} a_{i'} c_l$ appears, and again, as all the coefficients are positive
polynomials and $i' + l \not\equiv 0 \mod k$, it follows that $a_{i'}(1)c_l
(1) < \epsilon a(1)c(1)$. So $r a(1) c_l (a) /k < \epsilon a(1)c(1)$. This entails $r < \epsilon k$.

If we select $r = \epsilon k$, we obtain a contradiction. Hence at most
one of $a_i$, say $a_j$, has norm at least $\epsilon k a(1)$, so that the
rest of the terms have norm strictly less. Then $\| a - x^j a_j \| <
(k-1)k\epsilon a(1)$. Similarly, with $c$, we obtain $j'$ \st $\| c -
x^{j'}c_{j'} \| < k(k-1) \epsilon c(1)$. If $j + j \neq k,0$, we see that
$a_j(1)c_{j'}(1) $ contributes completely to the error (as previously), so
$\epsilon a(1) c(1) > a_j(1)c_{j'}(1) > a(1) (1 - k(k-1)\epsilon)c(1)(1 -
k(k-1)\epsilon)$. This will lead to a contradiction if $\epsilon < (1 -
(k^2 - k )\epsilon)^2$. Then $\epsilon < 1/2(k^2 - k)$ is sufficient.

Hence $j+ j' = 0$ or $k$, and we are done.
\qed

\noindent In the following, all the inequalities boil down to optimization on
intervals or rectangles, and it always turns out that the minima are
achieved on the boundary; because the functions being minimized are at
worst quadratic, the arithmetic is quite simple.

\Lem Lemma \sevthr. Suppose $a,a' \in A^+$, with $\max \brcs{ a\proj i(1) }= \eta
a(1)$ and $\max \brcs{a'\proj i(1) }= \mu a'(1)$. If $\mu, \eta > 1/2$,
then $\dist (aa',B^+) /aa'(1) \geq \eta + \mu - 2\eta \mu >
\max\brcs{1-\mu, 1-\eta} $; otherwise, $\dist (aa',B^+) > 1/2$.

\Pf Write $a = \sum x^i a\proj i$ and $a' = \sum x^i a'\proj i$, so that
(with $ j \neq 0$)
$$\eqalign{
(aa')\proj 0 &= \sum a\proj i a'\proj {-i} \cr
(aa')\proj j & = \sum a \proj i a'\proj {j-i}\cr
}$$
Hence
$$\eqalign{
\dist (aa',B) &= \sum_{j\neq 0}\sum_i a \proj i (1)a'\proj {j-i}(1)\cr
& = a(1)a'(1) - \sum a\proj i(1) a'\proj {-i}(1)\cr
}$$
We may assume that $a(1) = a'(1) = 1$ (dividing by the appropriate
scalars). Fix $0\leq \eta, \mu < 1$, and consider the following problem.

Minimize $1 - \sum_{i=0}^{k-1} X_i Y_i$ subject to $X_i, Y_i \geq 0$,
$\sum X_i = \sum Y_i = 1$, and $X_j \leq \eta$, $Y_j \leq \mu$ for all
$j$.

[Set $X_i = a\proj i (1)$ and $Y_i = a'\proj {-i} (1)$.] In the interior
of the domain, the only critical point (after replacing $X_0$ and $Y_0$ by
$1 -\sum_{i\neq 0} X_i$ and
$1 -\sum_{i\neq 0} Y_i$ respectively)---$X_i = Y_j = 1/k$---yields the
maximum. So the minimum must occur on the boundary.

Now assume $\eta, \mu > 1/2$. It easily follows that the minimum occurs
only when for some $i$, $X_i = \eta$ and $Y_i = \mu$ (the same $i$), and
then only when there exist some other $j \neq i$ \st $X_j = 1 - \eta$ and
$Y_j = 1-\mu$.

So the minimum value is $1 - \eta \mu - (1-\eta)(1-\mu) = \eta + \mu -
2\eta\mu$.

Dropping the condition $\mu, \eta > 1/2$, write $1 = r \eta + \psi$ and
$1 = s \mu + \phi$ where $r,s$ are positive integers, $ 0\leq \psi < \eta$
and $0 \leq \phi < \mu$ the minimum occurs when (up to relabelling), $X_1
= \eta$, $X_2 = \eta$, $\dots, X_r = \eta$, $X_{r+1} = \psi$, and $Y_1
=\mu = \dots = Y_s$, and $Y_{s+1} = \phi$.

Hence if $r = s$, the minimum value is $1 - r\eta\mu - (1-r\eta)(1-r\mu) =
r (\eta + \mu) - (r^2+ r) \eta\mu$. If $r\geq 2$, then the minimum value
(obtained by minimizing with the constraints $1/{(r+1)< \mu, \eta \leq
1/r}$) larger than $1-1/r$, which exceeds $1/2$.

If we assume $1/2 < \eta,\mu \leq 1$, then
$$
\eta + \mu - 2\eta\mu > \max\brcs{1-\eta, 1-\mu}.
$$

If $r < s$, the mimimum value is $1 - r\eta\mu > 1 - r/(r+1)(s+1) > 1/2$.
Similarly, if $s > r$, the minimum value is at least $1/2$.
\qed

Next, we deal with sums of products and elimination of small terms.
Suppose $a$ has two components, $a^{(0)}$ and $a^{(j)}$ \st each has mass
at least some number $a(1)\delta$. We obtain an estimate for $\dist
(aa',B^+)$ when $a', a \in A^+$. The coefficient of $x^t$ in $aa'$ is $\sum
a^{(i)}a'{}^{(t-j)}$; suppose $t \neq 0$. Then the error is at least
(summing over all $t \neq 0$) $a^{(0)}(1) \sum a'{}^{(t)}(1) + a^{(j)}(1)
\sum_{u\neq -j} a'{}^{(u)}(1) $. In particular, ${a'}^{(0)}(1)$ appears in
the second summand. Hence $\dist(aa',B^+) \geq \delta a(1)a'(1)$. Thus we
have the following.

\Lem Lemma \sevfou. Suppose $a,a' \in A^+$, and there exists $\delta > 0$ \st at
least two components of $a$ have mass at least $\delta a(1)$. Then $\dist
(aa',B^+) \geq \delta aa'(1)$.

\noindent {\it Proof of Proposition \sixthr.} For $a \in A^+$, define $\eta (a) = (a(1) -
\max_i\brcs{a^{(i)}}(1))/a(1)$. If $\eta (a) < 1/2$, then the $i_0$ for
which $a^{(i_0)}(1)$ is maximal is unique, and thus we can define $\pi(a)$
to be $x^{i_0} a^{(i_0)}$. It is helpful to keep track of the index $i_0$,
so we define the affiliated map $\Pi (a) = (a^{(i_0)}, i_0)$ with values
in $B^+ \times \Z_k$.
The error resulting from replacing $a$ by $\pi(a)$ is $\eta(a) a(1)$.

If the second largest of $\brcs{a^{(i)}}$ occurs at $i_1 \neq i_0$, and
$a^{(i_1)} = \delta$, then the previous result yields $\dist (aa',B^+)
\geq \delta aa'(1)$. However, the second largest is at least as large as
$(1 - \eta(a))/(k-1)$. Hence we deduce that in general, $d(aa',B^+) \geq
(1- \eta(a))aa'(1)/(k-1)$.

\comment
Roughly speaking, If $(1-\eta(a)/(k-1) > \epsilon$ (that is, $1 - \eta(a)
> (k-1)\epsilon$), then deleting $a$ actually improves the approximation!
We use this (with $\sqrt{\epsilon}$ replacing $\epsilon$) to show we can
delete all the bad matrix entries all at once, at a cost of
$\sqrt{\epsilon}$.
\endcomment

Fix a column of $V$, and pick one of its entries $a$ \st $(1-\eta(a)) >
(k-1)\sqrt{\epsilon}$. Suppose that $a$ is the $i$th entry. Let $c_t$ vary
over the entries of the $j$th column of $W$. Then the sum $\sum a c_j$ appears as part of
the sum of the entries of the $i$th column of $W$. Now sum this over all
the possible $a$ in the $i$th column (with the same condition on $\eta$),
obtaining $e:=\sum a_{s}c_j $, where $a_s$ varies over the entries in the
$i$th column of $V$ \st $(1-\eta(a_s)) > (k-1)\sqrt{\epsilon}$. Then
$\dist (e,B^+) \geq \sqrt {\epsilon} e(1) \sum c_j(1)$.

If we knew that $\sum c_j (1) = 1$ for all $j$, we would deduce $e(1) <
\sqrt {\epsilon}$. Now do this simultaneously for all the columns of $V$;
since the norm is the maximum column sum we deduce that simultaneously
replacing all the entries satisfying the condition on their $\eta$ by zero
increases the norm of $V$ by at most $\sqrt{\epsilon}$.

Hence at a cost of increasing $\epsilon$ to $\sqrt{\epsilon} + \epsilon$,
we may suppose that all the nonzero entries of $V$ satisfy $(1- \eta(a))
\leq (k-1)\sqrt{\epsilon}$. The next stage is to replace each entry by
$\pi(a)$, the projection onto the $X^0$-component (out of the $X^0, \dots , X^{k-1}$ components). The error along each column is just the sum of $a_{i}(1) (1-
\eta(a_i))$, which is bounded above by $(k-1) \sqrt{\epsilon}$ (since the
column sum of $a_i(1)$ is $1$). Hence at a cost of an additional
$(k-1)\sqrt{\epsilon}$, we can assume that the entries of $V$ are all of
the form $x^{f(i,j)}b_{ij}$ where $\Arrow f; \N \times \N .
{\brcs{0,1,2,\dots, k-1}}$. (The cumulative error on $V$ is now
$k\sqrt{\epsilon} + \epsilon$.)

Since $WV $ is close to a matrix $N$ with entries in $B^+$, we can of
course do the same thing to $W$. This allows us to assume (with extra
error introduced by changing $W$, $k \sqrt{\epsilon} + \epsilon$). Let
$\kappa = \max\brcs{\| V'W' - M \| , \| W'V', N\|}$, so that $\kappa < 2k
\sqrt{\epsilon} + 2 \epsilon$. Write the entries of $V$, $v_{ij} =
x^{f(i,j)}b_{ij}$, $w_{jl} = x^{g(j,l)})b'_{jl}$ where the $b$s all lie in
$B^+$ and the ranges of $f$ and $g$ are in $\brcs{0,1,2,\dots, k-1}$.

Now assume that $W_i, V_i$ with appropriate summable errors $\| W_iV_i -
M_i \| = \epsilon(i)$, etc are defined, and we may assume (by something
above) that $\sum \sqrt{\epsilon(i)} < \infty$. By the previous
construction, we can assume the entries of $W_i$ are all of the form
$(W_i)_{jl} =x^{ g_i(j,l)} b_{ijl}$ and $(v_i)_{jl} =x^{ f_i(j,l)}
b'_{ijl}$ where $f_i,g_i$ have ranges in $\brcs{0,1,2, \dots, k-1}$. We
may replace the $f_i$ by $f'_i$ where the range of each $f'_i$ is in
$\brcs{0,-1,-2, \dots, -(k-1)}$.

At this stage, we have a trick. Replace all the entries of all the
$V_i$ and $W_i$ by removing the $x^{f_i(j,l)}$ and $x^{g_i(j,l)}$ terms
(so the entries of the new $W_i$ and $V_i$ are respectively $b_{ijl}$ and
$b'_{ijl}$). Now it is almost trivial that that the resulting errors in
approximating $M_i$, respectively $N_i$ by products $M_i -W'_i V'_i$ and
$N_i - V'_{i+1} W'_i$ do not increase! The one-sided version (where we do not assume the existence of the $N_i$, but merely define $N_i:= V'_{i+1} W'_i$) of this of course applies, and will the AT($n$) results in part (b).

To see this, we look at an individual term in an individual entry; if
$f_i(j,l) + g_i (l,m) \neq 0$ (because of our assumption on the values for
the $f_i$ and $g_j$s, we cannot obtain $f_i(j,l) + g_i (l,m)= \pm k$), then
the corresponding term yields an error (in the original factorization) of
$b_{ijl}b_{lm}(1)$ to that term (since it is part of the distance to $B$);
hence when we remove the $f_i$s and $g$s, the error will be at most that,
and since the errors in the distance to $B$ are additive, we find that the
error resulting from this contribution will be at most that from the
original. On the other hand, if $f_i(j,l) + g_i(j,l) = 0$, the corresponding
term will be exactly the same as that in the original factorization. So
the outcome is that $\| W_i V_i - W'_i V'_i\| < 2k(\sqrt{\epsilon(i)} + 2 \epsilon(i)$ and similarly (if necessary)
for $\| V_iW_{i+1} - V'_i W'_{i+1}\| $.
\qed

\comment

\thmxx do the more general results: applying $Q$ or $xP$ yields $A$-iso,
and $A$-module iso between $B$ modules implies $B$-module iso. Clarify notion of iso over $A$ vs $B$.

\Pf Suppose that $M$ is an $r \times s$ matrix with entries from $B^+$, and there exists an approximate factorization (just how approximate will be dealt with later), $\| M - VW \| < \epsilon$ where $V = (v_i)^T \in (A^{r \times 1})^+$, $W = (w_j) \in (A^{1 \times s})^+$, and the norm is the maximum of the column sums of $l^1$-norms; for example, $\| V \| = \sum \| v_i\|$, the latter as elements of $A \subset l^1(\Z)$.
We also assume that the column sums of $M(1)$ are all $1$, so that $(1\ 1 \ \dots \ 1)$ is a left eigenvector for $M(1)$, with eigenvalue $1$. Then with error at most $\epsilon$ and working in finite-dimensional real $l^1$ spaces, $(1\ 1 \ \dots \ 1) V(1)W(1) \sim (1\ 1 \ \dots \ 1)$, whence by dividing $V$ by the positive real number $(1\ 1 \ \dots \ 1) V(1)$ and multiplying $W$ by the same number, we can assume at the outset that $W(1) \sim (1\ 1 \ \dots \ 1) $, and $\sum v_i (1) \sim 1$.

Recall that for $a$ in $A^+$ or $B^+$, $\| a\| = a(1)$ (since $x \mapsto 1$ is equivalent to $X \mapsto 1$ and the image of $x$ is positive, $a(1)$ is unambiguous when $a \in B^+$).

With the new choice of $V$ and $W$, we may assume that $W(1) = (1\ 1 \ \dots \ 1)$
and $\sum_i v_i(1) = 1$, and now $\| M- VW \| < 3\epsilon:= \kappa$. Write $M = (m_{ij})$. We have, for each $j$, $\sum_i \| m_{ij} - v_i w_j \| <\kappa$. Now introduce another small number, $\eta$, which will be determined later, except that we have to impose the mild condition, $(k^2 - k)\eta < \slfrac 12$.

Let $S_j = \Set{i \in {1,2,\dots,r}}{\| m_{ij} - v_i w_j \| \geq \eta v_i(1)}$. If $i$ in ${1,2,\dots,r}$ is not in $S_j$, then by the preceding lemma, there exists $t \equiv t(i,j)$ together with $v_i'$ and $w_j'$ in $B^+$ \st both $ \| v_i - x^t v_i' \| <(k^2 - k ) \eta v_i(1)$ and $\| w_j - x^{-t} w_j\| < (k^2 - k ) \eta w_j(1) = (k^2 - k ) \eta$. On the other hand,
$$
\kappa> \sum_{i \in S_j} \| m_{ij} - v_i w_j\| \geq \eta \sum_{S_j} v_i (1).
$$
Thus $\sum_{S_j} v_i (1) < \kappa/\eta$. In particular, $\sum_{i \in S_j} \| m_{ij} (1) \|< 2\kappa/\eta$. Replacing each $v_i$ by zero for all $i$ in $\cup_j S_j$ results in an error of at most $2\kappa /\eta$ in the approximate factorization.

For the remaining $i$, we wish to show that $t (i,j)$ (the exponent) is independent of the choice of $i$ and $j$. This is an easy consequence of the method of the previous lemma.

Explicitly, suppose $a = \sum_{l = 0}^{k-1} x^l a_l$ with $a_l$ in $B^+$, and there exists $t' \in \brcs{0,1,2,\dots,k-1}$ together with $a'$ in $B^+$ \st $\|x^{t'} a' - a\| < \alpha a(1)$. Then
$$
\| a' - a_{t'} \| + \sum_{l \neq t'} a_l (1) = \|x^{t'} a' - a\| < \alpha a(1).
$$
The second summand is $a(1) - a_{t'} (1)$, so we deduce that $\| a' - a_{t'} \| < a_{t'} (1) - a(1) (1- \alpha)$. In particular, $a_{t'} (1) > a(1) (1-\alpha)$. As $sum_{l} a_l (1) =1 $, for all other values of $t \neq t'$, we must have $\sum_{l \neq t'} a_l (1) < \alpha a(1)$. If $\alpha < \slfrac12$, uniqueness of the $t'$ is immediate.

Hence, provided $(k^2 - k) \eta < \slfrac 12$ (which we assumed at the outset), the choice of $-t$ is unique (this time in $\brcs{0,-1,-2, \dots , -(k-1)}$) for each $w_j$ and therefore for all the remaining $i$.

Set $v_i' = 0$ if $i \in \cup S_j$; the remaining values for $v_i'$ and $w_j'$ are already determined, as has $t$; this yields $V' = (v_i')^T$ in $(B^{r \times 1})^+$ and $W' = (w_j')$ in $(B^{1 \times s})^+$. Then
$$\eqalign{
\| x^{-t}V - V' \|& \leq (k^2 - k) \eta \sum_{i \not\in \cup_j S_j} v_i(1) + \frac{\kappa}{\eta} \leq (k^2 - k) \eta + \frac{\kappa}{\eta} \cr
\| x^{t} W - W' \| & \leq \max_j \brcs{(k^2 - k) \eta w_j(1)} = (k^2 - k) \eta; \qquad \text{thus} \cr
\| M - V'W' \| & \leq \| M- VW\| + \| V'\| \(\| x^t W - W' \) + \| W' \| \( \| x^{-t}V - V'\|\) \cr
& \leq \kappa + \(1 + (k^2 - k)\eta + \frac {\kappa}{\eta}\) (k^2 - k) \eta + \(1+ (k^2 - k) \eta \)\( (k^2 - k) \eta + \frac{\kappa}{\eta}\). \cr
}$$
If we insist that $\epsilon < (k^2 -k )/12$, then $\kappa < (k^2 -k)/4$, and now we set $\eta = \sqrt{\kappa/(k^2 -k)}$. Then $\eta < 1/2(k^2 -k)$ (as was required), and the final error term goes to zero as $\epsilon \to 0$.

The first result is now immediate. The second has (after telescoping) $W_{u+1}V_u \sim p_u$, and it easily follows that the choice of $t$ is even independent of the indexing, so $(W'_t)$ and $(V'_t)$ implement a $B$-module isomorphism.
\qed
\endcomment
\noindent {\it Remark} In other words, if two $B$- (or $l^1(k\Z)$-) dimension spaces become isomorphic on tensoring over $B$ with $A$ (or $l^1 (\Z)$), then they were isomorphic to begin with.

As a sample computation (the basis for the general odometer result), with an extra point, consider the sequence of matrices generated by $T^2$, where $T$ is the $3$-odometer. Here $T$ is represented by the sequence of polynomials $(p_j:= (1+ x^{3^j} + x^{2\cdot 3^j})/3)$. Then $\Cal B(p_j) = \slfrac13 \(\smallmatrix 1 + X^{3^j} & X^{(3^{j}+1)/2} \\ X^{(3^{j}-1)/2} & 1 + X^{3^j} \endsmallmatrix\)$. Products of these matrices can be calculated directly, but instead, let us consider the corresponding matrices over $A$ (instead of $B$).

Replacing $X$ by $x^2$, and applying conjugation with $\Delta$, we obtain a sequence $\(M_j = \slfrac13 \(\smallmatrix 1 + (x^{3^j})^2 & x^{3^{j}} \\ x^{3^{j}} & 1 + (x^{3^j})^2 \endsmallmatrix\)\)$. This is an ergodic circulant sequence, and the eigenvalues of $M_j$ are $(1\pm x^{3^j} + x^{2\cdot 3^j})/3$. It easily follows that $(M_j)$ is not hollow (this is always true when the sequence is derived from ergodic $T^2$ in this fashion). We can see that $(M_j)$ is AT (without referring to $\Cal B(p_j)$) by obsering that $\prod_{j=0}^{N-1} M_j = 3^{-N}\(\smallmatrix P_0 & P_1 \\ P_1 & P_0\\\endsmallmatrix\)$, where $P_0 = \sum_{\text{even }i < 3^N} x^i$ and $P_1 = \sum_{\text{odd }i < 3^N} x^i$, and since $\|x^{\pm1}P_1 - P_0 \| = 2$, we can conclude that $(M_j)$ is AT, and now the obvious approximate factorization yields that $(M_j) \iso ((1 + x^{2\cdot 3^j} + x^{4\cdot 3^j})/3)$, exactly as we would have obtained from $(\Cal B(p_j)) \iso ((1 + X^{3^j} + X^{2\cdot 3^j})/3)$ by replacing $X$ by $x^2$.

Note however, that $(M_j)$ is not isomorphic to the 3-odometer (we can apply mass cancellation invariants). Its dimension group corresponds to the odometer with supernatural number $2\cdot 3^{\infty}$, which although Kakutani equivalent to the $3$-odometer, is not conjugate to it. (We must remember that isomorphisms of $(\Cal B(p_j))$ are implemented over $B = \R[X^{\pm1}]$, while isomorphisms of $(M_j)$ are implemented over $A$.)

As a slightly different example, consider what happens when $T$ is represented by $(p_j = (1+ x^{g(j)})/2)$ as in Proposition \twoeig. Then
$$
\Cal B(p_j) = \cases \frac 12 \((1+ X^{g(j)/2})\I\) & \text{if $g(j)$ is even}\\
\frac 12 \(\I + X^{(g(j)-1)/2}\(\smallmatrix 0 & X\\ 1 & 0\\ \endsmallmatrix\)\) & \text{if $g(j)$ is odd.}\\
\endcases
$$
It follows that $T^2$ is ergodic if and only if infinitely many $g(j)$ are odd.

If, for example, {\it all\/} the $g(j)$ are odd, then $\Delta \Cal B (p_j)\Delta^{-1} = M_j = \slfrac12\(\smallmatrix 1 & x^{g(j)}\\ x^{g(j)} & 1\\ \endsmallmatrix\)$, so that if we assume the conditions on $g(j)$ of \twoeig, then $(M(j))$ is AT, and thus so is $(\Cal B(p_j))$, the latter by Proposition \sixthr.

If instead $g(j)$ is the $j$th Fibonacci number (counting $2$ as the third one), then $g(3j)$ is even; when we telescope in threes, we obtain
$$
\(M_{3j}M_{3j+1}M_{3j+2} = \frac18 (1+ x^{g(3j)})\(\smallmatrix 1 & x^{g(3j+1)}\\ x^{g(3j+1)} & 1\\ \endsmallmatrix\)\(\smallmatrix 1 & x^{g(3j+2)}\\ x^{g(3j+2)} & 1\\ \endsmallmatrix\)\).
$$
This is not hollow (in contrast to the corresponding sequence considered in Lemma \twosix), but is AT (this is true much more generally). Therefore $T^2$ is AT.

We can ask whether $T$ and $T^2$ are conjugate (when the latter is ergodic), that is, whether their corresponding dimension spaces are isomorphic. This is trickier, especially since it happens to be true for odometers, and we expect it won't be for generic AT (and other) actions.

If in the preceding example, we set $g(j) = 5^j$ (so none of the previous results apply), then the corresponding $\(M_j = \slfrac12\(\smallmatrix 1 & x^{5^j}\\ x^{5^j} & 1\\ \endsmallmatrix\)\)$ is obviously ergodic, but is known [GH] {\it not\/} to be AT. Therefore, by \sixthr, $T^2$ is not AT, but is ergodic.

Now we have a simple but useful result on non-hollowness of sequences of the form $(p_i(xP))$ that arise from conjugation of $(\Cal B(p_i))$ by $\Delta$.

If $S$ and $T$ are ergodic transformations represented respectively by $(M_j)$ and $(N_j)$, then we define the transformation $S \otimes T$ to be the transformation determined by the dimension space obtained from $(M_j \otimes N_j)$. When we assume (as we have throughout) that there is a unique invariant measure for each of $S$ and $T$, then it easily follows that the corresponding trace on the dimension space of $(M_j \otimes N_j)$ is ergodic; hence, at least when there is an invariant measure for each, $S\otimes T$ is ergodic. Unfortunately, there does not seem to be a dynamical characterization of $S\otimes T$, even when $S = T$. In the situation of topological dynamics, at least for Vershik's adic transformations, $S\otimes T$ has been given a dynamical meaning [BH], but it is not entirely satisfactory.

\Lem Lemma \sevfiv. Suppose that the AT transformation $T$ is represented by the sequence of polynomials $(p_i)$, and suppose that $T^k$ is ergodic. Let $P$ denote the standard cyclic permutation matrix of size $k$. Then for all positive integers $l$, the ergodic circulant sequence (of size $k$ matrices) $((p_i(xP))^l)$ is not hollow.

\Pf Since $T^k$ is ergodic, the sequence of matrices over $B^+$, $(\Cal B(p_i))$ is ergodic. Conjugating each of the terms $\Cal B(p_i) = p_i (Q)$ with $\Delta$ and replacing $X$ by $x^k$, we obtain a new sequence of matrices from $A$ given by $(p_i(xP))$. Ergodicity of the former sequence implies that of the latter (since the criteria for ergodicity just depend on the real matrices obtained by assigning $x$ and $X$ to $1$).

With $H = \Z_n$, the dual group runs over the $k$ characters, $\chi_j$, sending the generator to $\xi^j$ where $\xi = \exp 2\pi \sqrt{-1}/k$; index the corresponding eigenvalue functions $\lambda_j$, $j \in \Z_k$. Now $\lambda_{j}(p(xP)) = p(x\xi^j)$, and thus $\lambda_j (\prod p_i(xP)) = p(x\xi^j)$ where $p = \prod p_i$. If we write $p = \sum c_t x^t$ with $\sum c_t =1$ and $c_t \geq 0$, then $p(x\xi^j) = \sum c_t \xi^{tj} x^t$, so that $ \| p(x\xi^j) \| = \sum |c_t \xi^{tj}| = \sum c_t = 1$. Thus $\lambda_j \(p(xP)\) = \prod \lambda_j\(p_i (xP)\)$ has norm one. Hence $(p_i (xP))$ is not hollow.

Now $T^k \otimes T^k \otimes \dots \otimes T^k$ (with $l$ factors) is ergodic (tensor products of ergodic transformations are ergodic, unlike the situation for cartesian products), and is the $k$th power of $T \otimes T \otimes \dots \otimes T$, which in turn is represented by $(p_i^l)$; hence the result of the previous paragraph applies to $(p_i^l (xP) = (p_i(xP))^l)$.\qed

Combined with earlier results, this says that $(p_i (xP))$ arising in this construction is not isomorphic to the original sequence $(p_i)$ (although it can be isomorphic to $(p_i(x^k))$, as occurs for odometers).

The following is an easy consequence of earlier results. However, its significance is reduced by the fact that we really don't know what $T\otimes T$ is dynamically.

\Lem Proposition \sevsix. Suppose $T$ is AT and for a positive integer $k$, $T^k$ is ergodic. Then $T^k \otimes T^k$ is AT.


\long\def\Rf[#1] #2, #3. #4\par%
{\vskip 2pt \itemitem{[#1]} #2, {\it #3,} #4\par\vskip2pt}

\SecT References

\Rf [Al] EM Alfsen, 	
Compact convex sets and boundary integrals. Springer--Verlag 1971, ix + 210\, pp.

\Rf [AE] L Asimow \& AJ Ellis, Convexity theory and its applications in functional analysis. London Mathematical Society Monographs 16 (1980).

\Rf [BH] BM Baker and DE Handelman, Positive polynomials and time dependent integer-valued
random variables. Canad J Math 44 (1992) 3--41.

\Rf [CW] A Connes \& EJ Woods, Hyperfinite von Neumann algebras and Poisson boundaries of time-dependent random walks, Pac J Math 137 (1989) 225--243.

\Rf [DQ] A Dooley \& A Quas, Approximate transitivity for zero-entropy systems. Erg Thy \& Dyn Sys 25 (2005) 443--453.


\Rf [EG] GA Elliott \& T Giordano, Amenable actions of discrete groups. Erg Thy \& Dyn Sys 13 (1993) 289--318.

\Rf [GH] T Giordano and D Handelman, Matrix-valued random walks and variations on property AT. MŸnster J Math 1 (2008) 15--72.

\Rf [GoH] KR Goodearl and DE Handelman, Metric completions of partially ordered abelian
groups. Indiana Univ Math J 29 (1980)  861Ð-895.

\Rf [H]   D Handelman, Isomorphisms and non-isomorphisms of AT actions. J d'analyse mathŽmatique, 108 (2009) 293--366.

\Rf [V] AM Vershik, Uniform algebraic approximation of shift and multiplication operators.
Dokl Akad Nauk SSSR (259) 526--529, 1981 (Russian).

{}

\vskip 10pt

Mathematics Department, University of Ottawa, Ottawa ON  K1N 6N5, Canada; dehsg\@uottawa.ca

\end